\documentclass[12pt, a4paper]{article}
\usepackage[left=1in,right=1in,top=1.0in,bottom=1in]{geometry}

\usepackage{amsmath,amssymb,amsthm,amsfonts,afterpage,hyperref}
\usepackage{amscd,subeqnarray}

\usepackage{graphicx}
\usepackage{caption}
\usepackage{subcaption}
\usepackage{authblk}
\usepackage{tikz}

\usepackage{array}
\usepackage{wrapfig}
\usepackage{color}
\usepackage{ragged2e}
\usepackage{stmaryrd}
\usepackage{siunitx}
\usepackage{multirow}
\usepackage{xcolor,cancel}
\usepackage{boldfonts}
\usepackage{symbols}
\usepackage{footmisc}
\usepackage{float}
\usepackage{setspace}
\usepackage{bm}
\usepackage{booktabs}
\usepackage{tikz, xcolor}
\usepackage{soul}

\usepackage[T1]{fontenc}
\usepackage[utf8]{inputenc}

\usepackage[linesnumbered,ruled,vlined]{algorithm2e}
\SetKwInput{KwInput}{Input}                
\SetKwInput{KwOutput}{Output} 

\usepackage[sort, numbers]{natbib}
\setcitestyle{square}

\usetikzlibrary{shapes,arrows, positioning,calc}	
\usetikzlibrary{arrows,snakes,backgrounds}

\newcommand{\bfa}[1]{\boldsymbol{#1}}

\DeclareMathAlphabet{\mathpzc}{OT1}{pzc}{m}{it}

\newcommand{\bfu}{\boldsymbol{u}}

\newcommand{\bfx}{\boldsymbol{x}}

\newcommand{\bfg}{\boldsymbol{g}}	

\newtheorem{theorem}{Theorem}[section]

\newtheorem{remark}{Remark}

\newtheorem*{tdiscretization}{Time discretization}
\newtheorem*{sdiscretization}{Weak Formulation}
\newtheorem*{fdiscretization}{Full discretization}

\providecommand{\keywords}[1]
{
  \small	
  \textbf{\textit{Keywords---}} #1
}

\title{A Finite Element Model for Hydro-thermal Convective Flow in a Porous Medium: Effects of Hydraulic Resistivity and Thermal Diffusivity}

\author[1]{S. M. Mallikarjunaiah\thanks{m.muddamallappa@tamucc.edu}\thanks{corresponding author}}
\author[2]{Dambaru Bhatta\thanks{dambaru.bhatta@utrgv.edu}}

\affil[1]{Department of Mathematics \& Statistics,
Texas A\&M University - Corpus Christi, 
Corpus Christi, Texas 78412-5825, USA}
\affil[2]{School of Mathematical \& Statistical Sciences,
The University of Texas - Rio Grande Valley, 
Edinburg, Texas 78539, USA}
\date{}

\begin{document}

\maketitle
	    
\begin{abstract}
In this article, a finite element model is implemented to analyze hydro-thermal convective flow in a porous medium.  The mathematical model encompasses Darcy's law for incompressible fluid behavior, which is coupled with a convection-diffusion-type energy equation to characterize the temperature in the porous medium.  The current investigation presents an efficient, stable, and accurate finite element discretization for the hydro-thermal convective flow model. The well-posedness of the proposed discrete Galerkin finite element formulation is guaranteed due to the decoupling property and the linearity of the numerical method. Computational experiments confirm the optimal convergence rates for a manufactured solution.  Several numerical results are obtained for the variations of the hydraulic resistivity and thermal diffusivity. In the present study, the bottom wall is maintained at a constant higher hot temperature while side vertical walls are thermally insulated and the top wall is maintained at a constant cold temperature. Heat transfer rates at the heated bottom wall are presented in terms of local Nusselt number. A linear variation in hydraulic resistivity and a quadratic variation in thermal diffusivity show an increase in the heat transfer rate.

 \end{abstract}
	 
\noindent \keywords{Finite element method; Hydro-thermal convection, Rayleigh number, Taylor-Hood elements, Porous Medium}

\section{Introduction}
Natural convection in a fluid environment is a process of upward movement of particles due to temperature difference and density variation or buoyancy force. Such convection is possible because of the change in properties of the elemental fluidic environment, while forced convection arises due to some external force. There are many cases of natural convection, for instance, in the classical  \textit{Rayleigh-B\'enard} convection (RBC) \cite{getling1998rayleigh,koschmieder1993benard,walker1978convection}, a slab of fluid when heated from below and cooled from above,  experiences temperature difference which leads to buoyancy force. The Darcy-B\'enard convection (DBC) problem is the porous medium equivalent of RBC. Hence, the fluid particles at the bottom become less dense in contrast to the particles at the top such a set-up leads to instability, later developing circular movement of the fluid particles forming convection cells or rolls \cite{chandrasekhar1981,ostrach1988natural}. In recent years, the RBC type convection is of great interest due to its significance in a wide range of engineering applications in cooling devices, manufacturing of energy storage devices, heat exchangers, crystal growth in liquids, air movement in closed environments such as airplanes, air-conditioned buildings, etc. 

There have been extensive studies done on buoyancy-driven convective flows and the evolution of heat and mass in enclosures. In particular, the simplest setting of steady laminar flow in a porous medium with differentially heated side walls is widely studied. The problem of characterizing the stability of onset of convection in DBC has been addressed by many researchers including issues with the boundary conditions at the bottom and top plates, controlling the instability using externally imposed forces \cite{horton1945convection,lapwood1948convection,walker1977note,homsy1976convective,rudraiah1980effects,nield1975onset,siddheshwar2019darcy,rees2000onset,alam2012numerical,muratov2021finite,si2023consistent}. The DBC problem with vertical variation of porous medium and fluid parameters are reported in \cite{bhatta2017convective,bhatta2020thermal}. It was reported that the positive vertical variations in hydraulic resistivity or thermal diffusivity have stabilizing effects in contrast to the negative variations of the same functions. However, a full numerical study is not achieved for the model presented in \cite{bhatta2017convective,bhatta2020thermal}. In the current paper, we propose a stable numerical discretization of the model reported in  \cite{bhatta2017convective,bhatta2020thermal} using the finite element method.

There are numerous works about the development of stable numerical discretization of the mathematical model describing the natural thermal convection in a porous medium such as projection methods, explicit finite elements, boundary element method, radial basis function collocation method \cite{ccibik2011projection,manzari1999explicit,grigoriev2017numerical,hickox1981numerical,vsarler2000natural,kosec2008local}. The primary aim of this paper is to introduce the standard continuous Galerkin based finite elements discretization of the natural thermal convection in a Darcy porous medium.  A stable Taylor-Hood element is used to stabilize the discretization and control the spurious oscillations from the pressure term discretization. Another major issue is the decoupling of the temperature term from the momentum equation, however, the velocity-temperature pair suffers time lag but inner iterations were performed until a stable stopping criterion is achieved. Overall, our algorithm and the subsequent implementation have proved efficient for a manufactured solution.  

The article is organized as follows: in the section~\ref{math_formulation}, a mathematical model describing the hydro-thermal convection in a porous medium is presented. Subsequently, a stable numerical method based on the continuous Galerkin-type finite elements is proposed. In Section~4, a convergence study by using a manufactured solution is presented and comparisons are also made about the heat transfer rate with the available studies.  A detailed study of the variation in fluid resistivity and diffusivity is presented in the section-4. Finally, a short conclusion is presented in the last section of the paper.

\section{Mathematical Formulation}\label{math_formulation}

This paper introduces an efficient finite element discretization of the hydro-thermal convection model that describes the behavior of a fluid when it experiences buoyancy differences due to temperature variations in a porous medium.  We consider a two-dimensional laminar flow of a non-isothermal, natural convection of viscous incompressible fluid bounded between two horizontal plates. In our mathematical model, we assume that the Bouusinesq approximation is valid, i.e. the density differences are ignored except in the term multiplied with the gravity.  Since the inertial force doesn't play any role in the fluid motion, the flow field at every time is entirely characterized by the buoyancy difference at that particular time. Hence the history of the flow field doesn't characterize the evolution of the fluid particles.  The temperature changes are large enough to affect the flow field because the fluid transports the heat and simultaneously the temperature field is affected by the changes in the flow field. This kind of coupling of the fluid flow and heat transfer is a common phenomenon in heat exchangers, and chemical reactors and also in instances in which the ingredients are cooled such as in alloy solidification.  In the current section, we first present the mathematical model for the hydro-thermal convection and in the subsequent sections, we propose a stable numerical method to obtain an approximate solution to the model.

 Let $\Omega \subset \mathbb{R}^2$ be a bounded polygonal domain with a Lipschitz boundary $\partial \Omega$ and $[0, \; T]$  be the time-interval with $T$ being the final-simulation time. The governing system of the partial differential equation system consists of a continuity equation for  the conservation
of mass, Darcy's approximation  for the conservation of momentum, and heat equation
for the conservation of energy. The non-dimensional system can be written in the primitive variables as \cite{nield2006convection,bhatta2017convective,fowler1997mathematical,bhatta2020thermal}:
\begin{subequations}
\begin{align}
\chi(z) \; \bfu -  {Ra} \; \theta\,\widehat{k}+\nabla p &=0 \quad \mbox{in} \quad \Omega \times (0, \; T),  \label{model_1a} \\
 \nabla\cdot \bfu  &=0   \quad \mbox{in} \quad \Omega \times (0, \; T),       \label{model_1b}           \\
 \frac{\partial \theta}{\partial t}+ \bfu\cdot\nabla \theta - \nabla \cdot  \left( \zeta(z) \;  \nabla \theta \right) &=  0  \quad \mbox{in} \quad \Omega \times (0, \; T),                         \label{model_1c}
\end{align}
\end{subequations}
where $\bfu \colon \Omega \times [0, \; T] \to \mathbb{R}^2$, $p \colon \Omega \times [0, \; T]  \to \mathbb{R}$, and $\theta \colon \Omega \times [0, \; T]  \to \mathbb{R}$ denotes the unknown velocity, pressure, and temperature fields, respectively; additionally $t$ denotes the time. Let $\bfa{x} = \left(x, \; z \right)$ is any point in the domain $\Omega$.  Further, ${Ra}$ is the \textit{Rayleigh number} which is a dimensionless quantity and it is the ratio of heat transport due to convection induced by buoyancy change and of heat transport due to thermal diffusion. The nondimensionlization introduced for the Darcy equation \eqref{model_1a} by using the following definition of  \textit{Rayleigh number}: 
\begin{equation}\label{Ra_def}
{Ra} := \dfrac{ g \, \beta \, \rho \, \delta \theta \, L^3}{ \eta \, \kappa},
\end{equation}
in which $g$ is the strength of the gravity vector in the downward $z$-direction, $\beta$ is the thermal expansion coefficient, $\rho$ is the density of the fluid, $\delta \theta$ is  the temperature difference, $\eta$ is the dynamic viscosity, $L$ is the reference length, $\kappa$ is the thermal conductivity. The functions $\chi(z)$ and $\zeta(z)$ are the hydraulic resistivity and the thermal diffusivity, respectively. 
The above mathematical model need to be supplemented by the initial and boundary conditions:
\begin{subequations}
\begin{align}
\bfu(\bfx, \; t) &= {\bfu}_D(\bfx, \; t), \quad \theta(\bfx, \; t) = {\theta}_D(\bfx, \; t),  \forall \; \bfx \in \partial \Omega, \quad \mbox{and} \quad t >0, \\
\theta(\bfx, \; 0) &= \theta^0(\bfx), \quad \forall \; \bfx \in  \Omega,
\end{align}
\end{subequations}

Note that the pressure variable doesn't need to have any boundary values, because on its own it doesn't have equations and there is no evolutionary structure. To overcome this difficulty, a strategy which is very popular in the literature,  \textit{penalty method} \cite{basak2006effects,reddy2019introduction,bhargava2010finite}, introduces pressure via a constraint equation:
\begin{equation}
p = - \gamma \; \nabla \cdot \bfa{u},
\end{equation}
where the constant $\gamma$ needs to be large so that a consistent solution can be obtained by satisfying the incompressibility criteria. In our current approach, for the \textit{uniqueness} pressure variable we assume that the discretization space for pressure satisfies the zero-mean normalization. Hence, we do not need to worry about fixing the improvised parameter ($\gamma$) for the stable numerical solution. 

In the next section, a numerical method based on the finite elements is presented for the stable discretization of the aforementioned  Darcy-B\`enard model with variations in the hydraulic resistivity and thermal diffusivity. 
\section{Numerical Method}
Let $\Omega \subset \mathbb{R}^2$ be an open, smooth, convex polygonal and connected domain. The boundary $\partial \Omega$ is assumed to be  Lipschitz continuous and has a non-zero measure.  Let $X$ be a space of real-valued functions defined on $\Omega$ with the standard inner product and $L^{p}\left( \Omega, \; X    \right)$ is the space of all \textit{Lebegue measurable functions} defined on $\Omega$. 
Then $L^\infty(\Omega)$ be the vector space of essentially bounded measurable functions with the usual  \textit{essential supremum norm}. 
In particular, $L^2(\Omega)$ is the space of all square-integrable functions with inner product $(f,\; g)_2  = \int_{\Omega} f \, g \; d\bfa{x}$ and norm $\| f \|_2 : = \sqrt{(f, \; f)_2 }$. For any $b \in \mathbb{N}, b \geq 0$, we denote by $W^{b, \, p}(\Omega)$ for the space of functions in $L^p(\Omega)$ admitting weak derivatives up to $p^{th}$ order, i.e.
\begin{equation}
W^{b, \, p}(\Omega):= \left\{ u \in L^{p}\left( \Omega  \right) \, \colon \; D^{\alpha}u := \dfrac{\partial^{|\alpha| }u}{\partial x_{1}^{\alpha_1} \ldots \; \partial x_{n}^{\alpha_n} }  \in  L^{p}\left( \Omega  \right), \;  |\alpha| \leq b, \right\},
\end{equation}
where $\alpha = \left( \alpha_1, \ldots , \alpha_n \right)$ is a multi-index order $|\alpha| = k$, and $k \in \mathbb{N}$.
In particular, 
\begin{equation}
H^1(\Omega) := W^{1,2}(\Omega) \quad \mbox{and} \quad H_0^{1}(\Omega) = \overline{C_0^{\infty}(\Omega)}^{\| \cdot \|_{H^1}},
\end{equation}
where $C^{m}(\Omega), \; m \in \mathbb{N}_0$ denote the linear space of continuous functions on $\Omega$. 
We equip the $H^1(\Omega)$ with a norm $\| \cdot \|_{H^1(\Omega)}$ and defined as 
\begin{equation}
\| u \|_{H^1(\Omega)} := \left( \| u \|_2^2 + \| \nabla u \|_2^2 \right)^{1/2}.
\end{equation}
We define a space $L^2_0(\Omega)$ which consists of functions in $L^2(\Omega)$ whose mean value is zero:
\begin{equation}
L^2_0(\Omega):= \left\{ v \in L^{2}(\Omega) \colon \int_{\Omega} v \; d\boldsymbol{x} = 0\; \right\}.
\end{equation}
\subsection{Time Discretization} \label{time_dis}
First, let us introduce the time-discretization as follows:  given a positive integer $N$, and the superscript $j$ indicates the  time-step  number:  
\begin{equation}
\delta t = t^{j} - t^{j-1} = \dfrac{T}{N}, \quad j=1, \ldots n, 
\end{equation}
which is a uniform time-step length. Let 
\begin{equation}
0=t^0 < t^1 < \; \ldots, \; < t^n = T
\end{equation}
be a discretization of the time interval $[0, \; T]$. 
The discrete-time formulation of the mathematical model reads as follows:
\begin{tdiscretization}
Let $\mbox{Ra} > 0$, the initial-boundary conditions $\bfa{u}_D(\bfa{x}, \; t)$, $\theta_D(\bfa{x}, \; t)$, and  $\theta^0(\bfa{x})$,   end-time  $T >0$, $\delta t = T /M $ be all given. Find the solutions $\bfa{u}^{n+1}  \in \left( H_{\bfg}^{1}(\Omega, \; [0, \; T]) \right)^{2}, \; p^{n+1} \in H^{1}(\Omega, \; [0, \, T]), \; \theta^{n+1} \in H^{1}(\Omega, \; [0, \, T])$ for all $n =0, \, 1, \, \ldots , (M-1)$, such that: 
\begin{subequations}
\begin{align}
\nabla \cdot \boldsymbol{u}^{n+1} &=0,  \label{cont-eq1} \\
\chi(z) \, \boldsymbol{u}^{n+1}  + \nabla p^{n+1} &= \mbox{Ra} \; \theta^{n} \,\widehat{k},    \label{vel-eq1} \\
  \theta^{n+1}+ \, \delta t \;  \boldsymbol{u}^{n+1} \cdot \nabla \theta^{n+1} - \, \delta t \; \nabla \cdot \left( \zeta(z) \nabla \theta^{n+1} \right) &= \theta^{n}.  \label{temp-eq1}
\end{align}
\end{subequations}
\end{tdiscretization}
In the above formulation, the time-derivative in the temperature equation is discretized by using \textit{forward Euler method}, which is a first-order accurate scheme. The other popular and higher-order methods are  \textit{backward differentiation formula}  or  \textit{Crank-Nicolson method} with some extrapolation \cite{girault2012finite}.  


\subsection{Weak Formulation }
In this paper, our goal is to solve the model problem proposed in Section~\eqref{math_formulation} and its corresponding time-discretization proposed in the previous subsection~\eqref{time_dis}. We use the finite element method for the spatial discretization of the time-discrete equations \eqref{cont-eq1}-\eqref{temp-eq1}. Some important reasons for the choice of finite elements over other discretization techniques such as finite difference or finite volume methods are the adjustability of the method with respect to the mesh refinement, and approximation over irregular geometry \cite{HLMSM2021,lee2020nonlinear,ferguson2015numerical,yoon2022preferential,yoon2022finite}.  More importantly, finite elements are powerful in controlling both discretization and approximation error by means of \textit{a posteriori} estimates \cite{oden1991hp}.

Before we present a weak formulation of the problem, let us define the abstract functional spaces to approximate the primitive variables such as velocity, pressure, and temperature,
\begin{align*}
\boldsymbol{X}_{\boldsymbol{0}} &:= \left( H^{1}\left( \Omega   \right) \right)^{2} = \left\{ \boldsymbol{v} \colon \Omega \to \mathbb{R}^2 \; | \; \boldsymbol{v} \in \left( L^{2}\left( \Omega \right) \right) ^{2}, \; \nabla \boldsymbol{v} \in  \left( L^{2}\left( \Omega \right) \right) ^{2 \times 2}, \; \boldsymbol{v}|_{\partial \Omega} = \boldsymbol{0}  \right\}  \\
\boldsymbol{X}_{\bfg} &:= \left( H^{1}\left( \Omega   \right) \right)^{2} = \left\{ \boldsymbol{v} \colon \Omega \to \mathbb{R}^2 \; | \; \boldsymbol{v} \in \left( L^{2}\left( \Omega \right) \right) ^{2}, \; \nabla \boldsymbol{v} \in  \left( L^{2}\left( \Omega \right) \right) ^{2 \times 2}, \; \boldsymbol{v}|_{\partial \Omega} = \bfg  \right\}  \\
X  &:=H^{1}\left( \Omega   \right) = \left\{ \phi \colon \Omega \to \mathbb{R} \; | \: \phi \in L^{2}\left( \Omega \right), \; \nabla \phi \in   \left(  L^{2}\left( \Omega \right) \right)^2  \right\}   \\
X_{\theta}  &:=H^{1}\left( \Omega   \right) = \left\{ \varphi \colon \Omega \to \mathbb{R} \; | \: \varphi \in L^{2}\left( \Omega \right), \; \nabla \varphi\in   \left(  L^{2}\left( \Omega \right) \right)^2, \; \varphi|_{\partial \Omega} = \theta_D  \right\}   \\
X_{p} &:= L^{2}_{0}\left( \Omega   \right) = \left\{ q \colon \Omega \to \mathbb{R} \; | \: q \in L^{2}\left( \Omega \right)  : \; \int_{\Omega} q \; d\boldsymbol{x}= 0  \right\},
\end{align*}
Each of the above functional spaces is equipped with an inner product and a corresponding norm defined at the beginning of this section. \\

To obtain a semi-discrete weak formulation, we multiply the time-discrete equations \eqref{cont-eq1}-\eqref{temp-eq1} by the test functions from the above functional spaces and integrate them by parts. Then, the weak formulation reads:

\begin{sdiscretization}
Let $\mbox{Ra} > 0$ the initial-boundary conditions $\bfa{u}_D(\bfa{x}, \; t)$, $\theta_D(\bfa{x}, \; t)$, and  $\theta^0(\bfa{x})$, end-time $T >0$,  $\delta t = T / M $ be all given.  Find the solutions $\bfa{u}^{n+1} \in \bfa{X}_{\bfg}, \; p^{n+1} \in X_p, \; \theta^{n+1} \in X_{\theta}$ for all $n =0, \, 1, \, \ldots , (M-1)$, such that: 
\begin{subequations}
\begin{align}
- \,  \left( q, \; \nabla \cdot \boldsymbol{u}^{n+1}   \right) &=0,  \label{wf-1}\\
 \left( \chi(z)  \; \boldsymbol{v}, \;  \boldsymbol{u}^{n+1}   \right)  - \left( \nabla \cdot \boldsymbol{v}, \;  p^{n+1}   \right)  &=   \mbox{Ra} \; \left( \boldsymbol{v},  \; \theta^{n} \widehat{k} \right),   \label{wf-2} \\
  \; \left( \phi, \; \theta^{n+1} \right) + \delta t  \; \left( \phi, \;  \boldsymbol{u}^{n+1} \cdot \nabla \theta^{n+1} \right) + \delta t   \; \left( \nabla \phi, \;   \zeta(z)  \nabla \theta^{n+1} \right) &=  \;  \left( \phi, \; \theta^{n} \right)   \label{wf-3}
\end{align}
\end{subequations}
for all $\boldsymbol{v} \in \boldsymbol{X}_{\boldsymbol{0}}, \; q \in X_p, \; \phi \in X$.
\end{sdiscretization}

\begin{remark}
It is easy to see that if the solution triple  $\left(\boldsymbol{u}, \; p, \; \theta  \right)$ satisfies \eqref{model_1a}-\eqref{model_1c} also satisfies the above continuous weak-formulation \eqref{wf-1}-\eqref{wf-3}. The converse can also be proved by using proper regularity assumptions on the test function spaces. 
\end{remark}

\begin{remark}
The existence of a solution to \eqref{wf-1}-\eqref{wf-2} can be shown as follows:  if $\Omega$ is a bounded, connected domain, then   there exists a constant $ \widetilde{c}>0$ satisfying the \textit{inf-sup} condition or {Ladyzhenskaya-Babu$\check{s}$ka-Brezzi}  (LBB) compatibility condition \cite{ladyzhenskaya1969}  at the ``continuous level'': 
\begin{equation}
\inf_{q \in X_p} \; \sup_{\boldsymbol{w} \in \boldsymbol{X}_{\boldsymbol{0}}} \dfrac{\left( q, \; \nabla \cdot \boldsymbol{w} \right) }{ \| q \|_{X_p} \; \| \boldsymbol{w} \|_{\boldsymbol{X}_{\boldsymbol{0}}} } \geq \widetilde{c} > 0,  
\end{equation}
then one can follow the usual arguments of the Lax-Milgram lemma for the existence and uniqueness of the solution to the equations \eqref{wf-1}-\eqref{wf-2} \cite{girault2012finite}. 
\end{remark}

\subsection{Fully Discrete Scheme}
In this part, we lay out an idea of stable finite element formulation of the continuous weak formulation described in the previous section. 
We discretize the domain $\Omega$ as follows:   
let $\mathpzc{T}_{h}$ be a subdivision of $\Omega$ into \emph{quadrilateral} elements $K \in  \mathpzc{T}_{h}$ such that $\overline{\Omega} = \bigcup_{K \in \mathpzc{T}_{h}} \overline{K}$ and $K_{i} \cap K_{j} = \varnothing$ for $K_{i}, \; K_{j} \, \in \mathpzc{T}_{h}$, with $i \neq j$. The finite element mesh is chosen in such a way that the boundary of $\Omega$ is discretized as the union of faces from  $\mathpzc{T}_{h}$. We assume that the subdivision of $\Omega$ is conforming: i.e. the intersection of the closure of any two elements $K_{i}$ and $K_{j}$ is either empty along a face. We also assume that this discretization is shape-regular (or regular) and simplicial in the sense of Ciarlet \cite{ciarlet2002finite}. Let $h_K$ be the diameter of each finite element $K$, and $h := \max_{K \in \mathpzc{T}_{h}} h_K$. The finite element space considered in this work is  $\mathcal{Q}_{k}, \; k \in \mathbb{N}$ and it is defined as: 
\begin{equation}
\mathcal{Q}_{k} := \left\{ {\varphi} \in C(\Omega) \,  \colon \; {\varphi} |_{K} \in \mathbb{Q}_k, \; K \in \mathpzc{T}_{h}  \right\},
\end{equation}
where $\mathbb{Q}_k$ is the tensor-product of polynomials up to order $k$ on the reference cell $\hat{K} = [0, \, 1] \times [0, \, 1]$. Then the discrete spaces to approximate velocity, pressure, and temperature variables are defined as follows. 
\begin{align}
\boldsymbol{X}_{h, \, 0} &= \left[\mathcal{Q}_{k+1}\right]^2 \, \cap \, \boldsymbol{X}_{\boldsymbol{0}}, \quad  \boldsymbol{X}_{h, \; \bfg} = \left[\mathcal{Q}_{k+1}\right]^2 \, \cap \, \boldsymbol{X}_{\bfg} \\
X_h &= \mathcal{Q}_{k+1} \, \cap \, X, \quad\quad X_{h, \; \theta} = \mathcal{Q}_{k+1} \, \cap \, X_{\theta}.  \\
X_{h, \, p} &= \mathcal{Q}_k \, \cap \, X_p.
\end{align}
For the approximation of velocity-pressure equations, \eqref{wf-1} and \eqref{wf-2}, we need to use the Taylor-Hood pair, $\boldsymbol{X}_h  \times X_{h, \, p} $, satisfies the \textit{discrete inf-sup} or LBB compatibility condition \cite{ladyzhenskaya1969}
\begin{equation}\label{lbb2}
\inf_{q_{h} \in X_{h, \, p}} \; \sup_{\boldsymbol{w}_h \in  \boldsymbol{X}_{h}} \; \frac{\left( q_{h}, \; \nabla \cdot \boldsymbol{w}_{h}  \right)}{ \| q_{h} \|_X \, \| \nabla \boldsymbol{w}_{h} \|_{\boldsymbol{X}}} \geq c >0,
\end{equation}
where constant $c$ is independent of $h$. A more detailed discussion on the  LBB condition can be found in \cite{girault2012finite}. For Darcy's equations, there are several possible stabilization choices for the approximation spaces for velocity-pressure pair that are compatible with discrete LBB condition \eqref{lbb2}, a simple choice is to use elements of one-order higher for velocities than for pressures. Otherwise, the choice of finite element approximation spaces will produce spurious pressure values and apparent velocities \cite{sani1981A}. Finally, the fully discrete formulation is given by:

\begin{fdiscretization}
Let ${Ra} > 0$, the initial-boundary conditions $I_{h} \left( \bfa{u}_D \right)$, $I_h \left( \theta_D \right) $, and  $\widetilde{\theta}^0$, end-time  $T >0$,  $ \delta t = \tau / N$ be all given.  Find $\bfa{u}_h^{n+1} \in \bfa{X}_{h, \, \bfg}, \; p_h^{n+1} \in X_{h, \, p},  \; \theta_h^{n+1} \in X_{h, \, \theta}$, for $n =0, \, 1, \, \ldots , (M-1)$ such that: 
\begin{subequations}
\begin{align}
-  \left( q_h, \; \nabla \cdot \bfa{u}_h^{n+1}   \right) &=0,  \label{ff-1}\\
 \left( \chi(z_h) \; \bfa{v}_h, \;  \bfa{u}_h^{n+1}   \right)  - \left( \nabla \cdot  \bfa{v}_h, \;  p_h^{n+1}   \right)  &=   \mbox{Ra} \; \left( \bfa{v}_h,  \;  \theta_h^{n} \widehat{k} \right),   \label{ff-2} \\
  \left( \phi_h, \;  \theta_h^{n+1} \right) + \delta t  \; \left(  \phi_h, \;  \bfa{u}_h^{n+1} \cdot \nabla  \theta_h^{n+1} \right) +  \delta t  \; \left( \nabla  \phi_h, \;   \zeta(z_h) \nabla  \theta_h^{n+1} \right) &=      \;  \left(  \phi_h, \;  \theta_h^n \right).   \label{ff-3}
\end{align}
\end{subequations}
for all $\bfa{v}_h \in \bfa{X}_{h, \, 0}, \; q_h \in X_{h, \, p} \;  \phi_h \in X_h$.
\end{fdiscretization}
In the above formulation, the initial value for temperature term i.e. $\widetilde{\theta}^0$ is the $L^2$-projection of the initial temperature onto the finite element mesh. Also,  $I_{h}(\cdot)$  is the appropriate interpolation operator. Notice that in \eqref{ff-2} we have decoupled the temperature, therefore temperature has a time lag. For better numerical accuracy, the equations \eqref{ff-1}-\eqref{ff-3} could be solved iteratively if the computational resources permit. 

\begin{theorem}
For any $\delta t >0$, the solutions to the above fully discrete finite element algorithm are unconditionally stable.
\end{theorem}

\noindent \emph{Proof: } Because of the linearity and finite-dimensional nature of the above proposed fully discrete finite element algorithm, the uniqueness follows from the linearity of the algorithm, and hence the existence \cite{girault2012finite}.

\subsection{Solution Strategy for the Linear and Decoupled Equations}
 The main goal of this paper is to develop an efficient solver for the hydro-thermal convection model described in Section~\ref{math_formulation}. To that end, we first need to formulate convenient linear algebra tools for the solution of finite-dimensional discrete problem proposed in the equations \eqref{ff-1}-\eqref{ff-3}. For the actual computations, we need to obtain the corresponding matrix-vector problem in terms of finite element basis functions. To obtain such matrices we need to explicitly write all the state variables as a linear combination of the basis functions. The finite-dimensional functional spaces defined by the basis sets as
\begin{align}
\boldsymbol{X}_h &:= \mbox{span}\left\{  \boldsymbol{\phi}_i            \right\}_{i=1}^{n_{\boldsymbol{u}}}, \quad
X_{h,p} := \mbox{span}\left\{  {\psi}_j            \right\}_{j=1}^{n_p},  \\
X_{h} &:= \mbox{span}\left\{  {\varphi}_k            \right\}_{k=1}^{n_{\theta}}, 
\end{align}
where $n_{\boldsymbol{u}}$, $n_p$, and $n_{\theta}$ are, respectively, the number of degrees of freedom for velocity, pressure, and temperature variables. The associated discrete solution variables, $\bfa{u}_h^{n+1} \in \left( \mathbb{R}^{n_{\bfa{u}}}\right)^2$, $p_h^{n+1} \in \mathbb{R}^{n_p }$, and $\theta_h^{n+1} \in \mathbb{R}^{n_{\theta}}$ are expressed as a linear combination of basis functions on the finite element mesh as
\begin{align}
 \boldsymbol{u}_h^{n+1}(t, \; x, \; z) &:= \sum_{i=1}^{n_{\boldsymbol{u}}} \boldsymbol{u}_i(t) \, \boldsymbol{\phi}_i(x, \; z), \\
 p_h^{n+1}(t, \; x, \; z) & :=  \sum_{i=1}^{n_p} {p}_j(t) \, {\psi}_j(x, \; z), \\
 \theta_h^{n+1}(t, \; x, \; z) & :=  \sum_{i=1}^{n_{\theta}} {\theta}_j(t) \, {\varphi}_j(x, \; z).
\end{align}

\subsubsection{Solver for velocity-pressure equations}
Some stable numerical methods are developed for the discretization of the coupled Darcy flow problem. Examples of such methods are the mixed finite element methods \cite{masud2002stabilized,correa2008unconditionally}, discontinuous Galerkin  \cite{brezzi2005mixed}, Reduced Basis Function collocation method \cite{kosec2008local}, finite volume method \cite{el2018analysis}. Our approach is very simple and straightforward, we discretize the PDE system for velocity and pressure variables by using \textit{Taylor-Hood} finite element pair, and bi-quadratic basis functions for temperature on the quadrilateral grid. 

In this section, we describe the efficient direct solver that we have implemented for the continuity-momentum pair of PDEs. First,  we need to compute the entries of the corresponding block matrices. We first consider the discretized velocity-pressure equations and write these equations using the basis functions:
\begin{align}
 \sum_{j=1}^{n_{\bfa{u}}} \bfa{u}_{j}(t) \left( \chi(z) \; \bfa{\phi}_i, \; \bfa{\phi}_j \right)
- \sum_{k=1}^{n_p} p_{k}(t)  \left(  \psi_k, \;  \nabla \cdot  \bfa{\phi}_i \right) &= {Ra} \, \left( \bfa{\phi}_i \cdot \widehat{k}, \;  \theta_h^n  \right), \quad i = 1, \, 2, \, \ldots n_{\bfa{u}} \label{dis_sys1} \\
- \sum_{j=1}^{n_{\bfa{u}}} \bfa{u}_{j}(t) \left( \psi_l, \; \nabla \cdot \bfa{\phi}_j \right) &=0, \quad 
l = 1, \, 2, \, \ldots n_p. \label{dis_sys2}
\end{align}
The above set of linear system of size $(n_{\boldsymbol{u}} + n_p) \times (n_{\boldsymbol{u}} + n_p)$ needs to be solved at each time-step for the unknowns $\boldsymbol{u}_{j}(t)$ and $p_{j}(t)$ at their respective degree-of-freedoms on the finite element mesh.  At every time step the assembled linear system \eqref{dis_sys1}-\eqref{dis_sys2} by using \textit{Taylor-Hood} $\mathbb{Q}_2$-$\mathbb{Q}_1$  finite element pair result in a matrix system of the form
\begin{equation}\label{linsystem1}
{\mathcal{A}} \widehat{\boldsymbol{x}} = G,
\end{equation}
where $\widehat{\boldsymbol{x}} = \left(U, \; P \right)^{T}$, right-hand side vector $G =  \left(F, \; 0 \right)^{T}$ and the left-hand side block matrix 
\begin{equation}\label{bmatrix1}
\mathcal{A} = \begin{pmatrix}
{A} &  & &B^{T} \\
B & &  &O
\end{pmatrix},
\end{equation}
where the vectors $U$ and $P$ are the discrete versions that belong to the finite element degrees of freedom of the velocity $\boldsymbol{u}_h(t, \, x, \, y)$ and pressure $p_h(t, \, x, \, y)$ variables. The block matrices are given by
\begin{subequations} \label{ABF_Msystem}
\begin{align}
A_{ij} &=  \left( \chi(z) \; \bfa{\phi}_i, \; \bfa{\phi}_j\right), \quad 1 \leq i, \; j \leq n_{\bfa{u}},   \\
B_{ij} &= \left(  \psi_i, \;  \nabla \cdot \bfa{\phi}_j \right), \quad 1 \leq i \leq n_{p},  \;\; 1 \leq j \leq n_{\bfa{u}},   \\
F_{i} &= {Ra} \, \left( \bfa{\phi}_i \cdot \widehat{k}, \;  \theta_h^n  \right), \quad 1 \leq i \leq n_{\bfa{u}}. 
\end{align}
\end{subequations}
The sparse linear system \eqref{linsystem1} is known in the literature as \textit{saddle-point system} and a great deal of work can be found in the literature about developing efficient solvers for such a system. Our goal here is to solve the matrix system obtained in \eqref{linsystem1}. The block structure in \eqref{linsystem1} can be solved for pressure first and then for velocity. The following are the equations to be solved for pressure and velocity.
\begin{subequations}
\begin{align}
BA^{-1}B^T \, P &= B A^{-1} \, F \\
U &= A^{-1} F - A^{-1}B^T P
\end{align}
\end{subequations}
In the above equations, we need to find $A^{-1} $, which is the inverse of velocity block, and then a matrix-vector multiplication with $B^T$.  To compute $A^{-1}$, we use \textit{Unsymmetric MultiFrontal method}  (UMFPACK) \cite{davis2006direct,davis2004algorithm} direct solver.

\subsubsection{Linear Solver for Temperature Equation}
In contrast to the solver for the Darcy system of equations, the solver for the temperature is much simpler and straightforward. Any numerical solution to the temperature convection-diffusion equation is to deal with the convection part of the governing equation in addition to the diffusion part. In our numerical model described in equations \eqref{ff-1}, \eqref{ff-2}, and \eqref{ff-3},  the convection part with the known ``velocity'' term  (i.e. $\boldsymbol{u}_h^{n+1} \cdot \nabla  \theta_h^{n+1}$) is rather much easy to handle due to knowing the unique velocity at the $(n+1)^{th}$ time-step. 

First, we need to write the matrix representation of the temperature equation. To that end, let us consider the discretized temperature equation: 
\begin{align}\label{tEquation}
 \left( \varphi_h, \; \theta_h^{n+1} \right) + \delta t  \; \left(  \varphi_h, \;  \boldsymbol{u}_h^{n+1} \cdot \nabla \theta_h^{n+1} \right) +  \delta t  \; \left( \nabla  \varphi_h, \;   \zeta(z) \nabla  \theta_h^{n+1} \right) =   \left(  \varphi_h, \; \theta^n \right).
\end{align}
Then, write the temperature variable $\theta_h^{n+1}$ as  a linear combination of basis-function $\varphi$ as:
\begin{equation}\label{tShpaeFun}
\theta_h^{n+1}(t, \; x, \; z)  =  \sum_{j=1}^{n_{\theta}} {\theta}_j(t) \, {\varphi}_j(x, \; z).
\end{equation}
Using \eqref{tShpaeFun} in \eqref{tEquation} we obtain the matrix representation for the temperature variable as:
\begin{align}
  \left(\varphi_i, \; \varphi_j \right) + \delta t \; \left( \varphi_i, \; \boldsymbol{u}_h^{n+1} \cdot \nabla \varphi_j \right) + \delta t \; \left( \nabla \varphi_i, \; \zeta(z) \nabla \varphi_j \right) =  \left( \varphi_i, \; \theta^n \right).
\end{align}
The above formulation can be written in matrix-vector form as 
\begin{equation} \label{matrix_temp_equ}
\left(  \left[ M \right] + \delta t  \;  \left[ B \right] + \delta t  \;  \left[ S \right] \right) \left\{ \theta_i \right\} =   f_i \, ,
\end{equation}
in which, the entries of each matrix and right-hand side are of the form
\begin{align}
M_{ij} &= \left(\varphi_i, \; \varphi_j \right), \quad 1 \leq i, \; j \leq n_{\theta},  \\
B_{ij} &= \left( \varphi_i, \; \boldsymbol{u}_h^{n+1} \cdot \nabla \varphi_j \right),  \quad 1 \leq i, \; j \leq n_{\theta},\\
S_{ij} &= \left( \nabla \varphi_i, \; \zeta(z) \,  \nabla \varphi_j \right), \quad 1 \leq i, \; j \leq n_{\theta}, \\
f_i &=  \left( \varphi_i, \; \theta^n \right), \quad 1 \leq i \leq n_{\theta}. 
\end{align}
 The matrix for the temperature equation is the sum of three matrices,  of which one is a symmetric positive definite, another one is a skew-symmetric positive definite, and the other is a symmetric positive semidefinite matrix. Hence, the result is also positive definite but not symmetric. As a consequence, inverting this matrix with the Krylov-based iterative \textit{Generalized Minimal Residual Method} (GMRES) \cite{saad2003iterative} solver along with a simple preconditioner is inexpensive compared to inverting the Darcy-velocity matrix. 
 
\begin{remark}
In any thermal convection problem, one of the important postprocessing tools to analyze is the \textit{heat transfer coefficient} in terms of the local Nusselt number ($Nu$) on the lower plate (i.e. $z=0$), a quantity that characterizes the average vertical heat transport,  and it is defined by \cite{otto2011rayleigh,grigoriev2017numerical}:
\begin{equation}
Nu = \bfa{n} \cdot \nabla \theta_h(\bfa{x})\big|_{z=0},
\end{equation}
where $\bfa{n}$ is the normal to the plane. We define the average Nusselt number at the bottom plate  as
\begin{equation}\label{Nu_formula}
Nu_b = \int_{\Gamma_b} - \dfrac{\partial \theta_h }{\partial z} \; dx. 
\end{equation}
\end{remark}

Finally, we summarize all the steps involved in our approach for the numerical solution of the hydro-thermal convection model as below: 
\begin{itemize}
\item[1. ] \textit{Domain, boundary conditions, and computational mesh}: A  coarse mesh is created for a square domain and all the boundary parts are ``numbered'' to apply the boundary conditions later. Each boundary part is subjected to no-slip velocity conditions along with isothermal temperature conditions. Two vertical walls are assumed to be adiabatic and have no normal flow conditions. The run-time parameters such as \textit{Rayleigh number}, hydraulic resistivity, and thermal diffusivity are read into the program from a parameter header file. 

\item[2. ] \textit{Initial temperature field}: The initial temperature filed $\theta^0(\bfa{x}) \in L^{2} \left( \Omega \right)$ is projected on to the finite element mesh. 

\item[3. ] \textit{Setting up of degrees of freedom and initialization of matrices and vectors}: Start the time-step loop and set up the \textit{degrees of freedom} and \textit{sparsity pattern} using $\mathbb{Q}_2$ and $\mathbb{Q}_1$ test function spaces. Then, initialize the matrices and vectors by the appropriate order from knowing the approximate value of the degrees of freedom on the computational mesh. 

\item[4. ] \textit{Solution to Darcy system and temperature equation}: Assemble  the matrices $A$, $B$ and the right-hand side vector $F$ as defined in \eqref{ABF_Msystem}, and create a linear system $\mathcal{A}\widehat{\boldsymbol{x}} = G$, where $\mathcal{A}$ is a $2 \times 2$ block-matrix, $\widehat{\boldsymbol{x}}$ and $G$ are $2 \times 1$ block-vectors. 
\begin{itemize}
\item[(i)] \textit{Linear solver for the Darcy system}: First solve the velocity-pressure linear system of equations by a direct method.  At this step, we obtain $\boldsymbol{u}_h^{n+1}(\bfa{x})$ and $p_h^{n+1}(\bfa{x})$ by using the old time-step temperature value $\theta^n(\bfa{x})$.
\item[(ii)] \textit{Solver for the temperature equation}: By knowing the velocity $\boldsymbol{u}_h^{n+1}(\bfa{x})$, solve the temperature equation \eqref{matrix_temp_equ} using \textit{GMRES} solver with a standard preconditioner.

\item[(iii)] \textit{Velocity-temperature inner iterations}: Solve velocity-temperature equations again by knowing the updated velocity and temperature values at the same time-step. Check whether the system has attained the steady state by an appropriate stopping criteria and then proceed to the next time step.
\end{itemize}

\item[5. ] \textit{Stopping criteria}: To find out whether the flow has reached the steady state, the following  stopping criteria are checked at all time steps:
\begin{equation}
\sqrt{\sum_{i, \, j} |\theta_{i\, j}^{n+1} - \theta_{i\, j}^{n}|^2} < 10^{-6},
\end{equation}
where superscript denote the time-step number and the subscript indices $(i\, j)$ denote the spatial coordinates. In all our numerical simulations, for a given Rayleigh number and other parameters in the fluid model,  the steady-state criteria are achieved and then the simulations are turned off. 

\item[6. ] \textit{Post-processing}: At each time step, save the variables such as velocity, pressure, and temperature on the entire computational mesh. At the end of the time-step loop, compute the \textit{Nusselt number} at the bottom plate.   

\end{itemize}

In the next section, we present some numerical tests both for the steady and unsteady Darcy-B\'enard convection problem. 

\section{Numerical Results and Discussion}

In this section, we present several numerical tests with the ultimate goal of validating the effectiveness of the proposed discretization of the Darcy-B\'enard convection model. We also present some comparison of our simulation results with the available results from the literature.  Several tests including checks for convergence rates of the numerical solution against a ``manufactured solution'' will be performed. Next, we consider a benchmark problem from the literature, and last tests were performed for several \textit{Rayleigh numbers} and also for different functions denoting hydraulic resistivity and thermal diffusivity. In all our numerical simulations,  \textit{inf-sup} stable \textit{Taylor-Hood} finite elements $\left( \mathbb{Q}_2, \; \mathbb{Q}_1\right)$ are used for velocity-pressure variables and a bilinear  $\mathbb{Q}_2$ elements are used to approximate temperature variable. The entire code is written using \textsf{C++} programming language based on an open-source finite element library  \texttt{deal.II} \cite{dealII91} along with \textsf{Visit} and \textsf{Paraview} for visualization. 

In the next few subsections, we list out a few numerical tests and present the results based on the finite element computations.  

\subsection{Convergence study}
In the first experiment, we present a convergence study to validate both the discretization of the mathematical model and its corresponding implementation. We consider some external forces $\bfa{f}(\bfa{x})$ as a right-hand side function in \eqref{model_1a} and $g(\bfa{x})$ in the steady-state version of the equation \eqref{model_1c} to use our numerical method the code for a particular ``manufactured solution''. 
\begin{subequations}\label{sm}
\begin{align}
 \bfu -  \theta\,\widehat{k}+\nabla p &= \bfa{f}(\bfa{x}) \quad\quad\quad \mbox{in} \quad \Omega,  \label{sm_1a} \\
 \nabla\cdot \bfu  &=0   \quad\quad\quad\quad\;\; \mbox{in} \quad \Omega,       \label{sm_1b}           \\
 \bfu\cdot\nabla \theta - \nabla \cdot    \nabla \theta  &=  g(\bfa{x})  \quad \quad\quad\; \mbox{in} \quad \Omega,                         \label{sm_1c}
\end{align}
\end{subequations}
For the convergence study, we choose a unit square $\Omega = (0, \; 1)^2$ and a manufactured solution as:
\begin{align}
\bfa{u} &= \begin{pmatrix}
\sin \pi x \\
- \pi \, z \, \cos \pi x
\end{pmatrix}, \quad  p = \sin \pi x \; \cos \pi z,  \label{exact_u_p}\\
 \theta  &= \sin(\pi \, x + \pi \, z). \label{exact_theta}
\end{align}
The above velocity vector satisfies the continuity equation. Then the corresponding right-hand functions $\bfa{f}(\bfa{x})$ and $g(\bfa{x})$ are calculated by substituting the known analytical solutions into \eqref{sm}. The exact forms of these two right-hand side functions are given by:
\begin{subequations}
\begin{align}
\bfa{f}(\bfa{x}) &= \begin{pmatrix}
\sin \pi x  + \pi \, \cos \pi x \; \cos \pi z \\
- \pi \, z \,  \cos  \pi x - \pi \, \sin \pi x \; \sin \pi z - \sin (\pi \, x + \pi \, z)
\end{pmatrix}, \\
& \notag \\
 g(\bfa{x}) &= 2 \pi^2 \, \sin(\pi \, x+ \pi \, z) + \pi  \, \cos(\pi \, x+ \pi \, z) \; \left(  \sin \pi x - \pi \, z \; \cos \pi x      \right)   
\end{align}
\end{subequations}

A total of $ 6$ cycles of global refinements are done and the finite element solutions are obtained for each refinement. The $L^2$ errors in each variable are recorded. The following table presents the errors for the velocity, pressure, and temperature variables. For solving the momentum equation, we choose the exact temperature analytical solution as the initial value. After obtaining the velocity vector, we compute the temperature values on the computational mesh. The boundary conditions for both velocity and temperature variables are set using the exact analytical solutions.

\begin{table}[H]
\centering
\caption{Errors in $L^2$- norms for velocity ($\bfa{u}$), pressure ($p$), and temperature ($\theta$) with the manufactured solution and their convergence rates.} 
\label{tab:h-conv}
\begin{tabular}{|c|| l|l|| l|l|| l|l|}
\hline
\multirow{2}{*}{DoFs}  & \multicolumn{2}{c||}{ $ \| \bfa{u}_h - \bfa{u} \|^2_{2}$}    & \multicolumn{2}{c||}{ $\| {p}_h - p \|^2_{2}$ }    & \multicolumn{2}{c|}{$\| {\theta}_h - \theta \|^2_{2}$}                            \\ \cline{2-7} 
&  \multicolumn{1}{c|}{$L^2$ Error}  & \multicolumn{1}{c||}{Rate} & \multicolumn{1}{c|}{$L^2$ Error} &  \multicolumn{1}{c||}{Rate} &  \multicolumn{1}{c|}{$L^2$ Error}  & \multicolumn{1}{c|}{Rate} \\ \hline\hline
187    & 0.162432  & 0.0    & 0.0174858 & 0.0    & 0.00364343  & 0.0           \\ \hline
659    &  0.0468909  & 1.8611   &  0.00414257 & 2.05450  & 0.00035787    & 3.19074   \\ \hline
2967   & 0.0124804  & 1.93834  & 0.00102083 & 2.01445  & 3.83226e-05  & 3.05587             \\ \hline
9539  &  0.00320821 & 1.97234   & 0.000254295    & 2.00358  & 4.58035e-06  & 2.892531 \\ \hline
37507  & 0.000812492  & 1.9871  &   6.35168e-05  & 2.00089 & 5.69871e-07    & 2.83505  \\ \hline
148739  & 0.000204388  & 1.99380  &  1.58756e-05 & 2.000226 &  7.14425e-08    & 2.82429   \\ \hline
\end{tabular}
\end{table}

The convergence rates in $L^2$-norms for the discretization using Taylor-Hood  $\mathbb{Q}_2\mathbb{Q}_1$ elements for velocity and pressure, as well as for temperature   
presented in Table~\ref{tab:h-conv}. It is clear from the above table that the discretization that we have done for the steady Darcy-B\'enard convection model, the convergence rates are optimal \cite{girault2012finite}. 

\begin{figure}[H]
\centering
\begin{subfigure}{0.3\textwidth}
    \includegraphics[width=\textwidth]{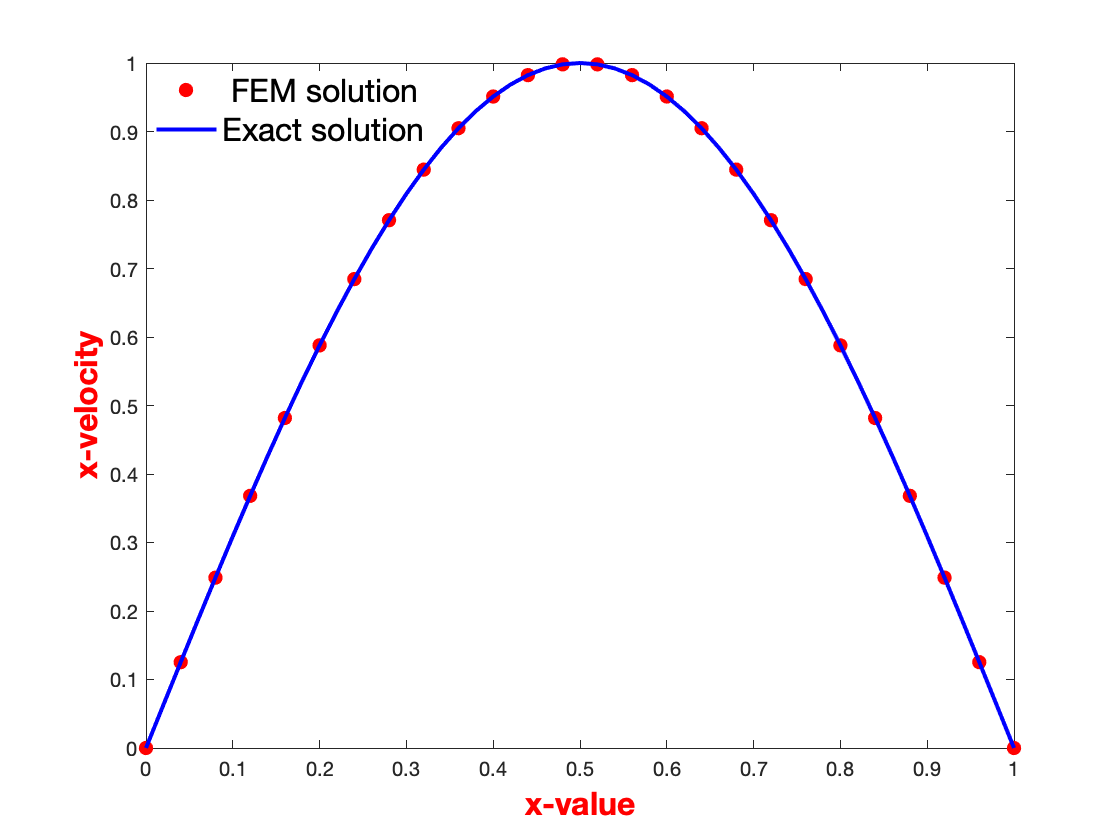}
    \caption{Plot of ${u_h}$}
    \label{fig:first}
\end{subfigure}
\quad 
\begin{subfigure}{0.3\textwidth}
    \includegraphics[width=\textwidth]{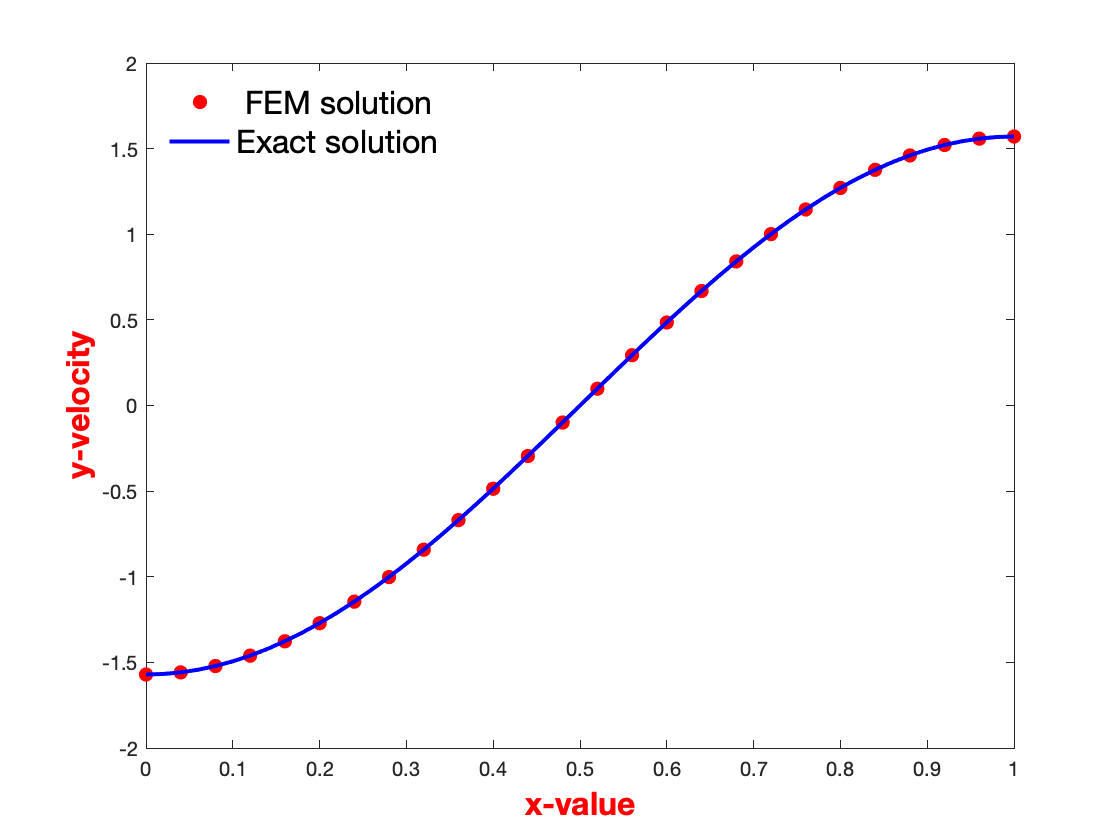}
    \caption{Plot of ${w_h}$}
    \label{fig:second}
\end{subfigure}
\quad 
\begin{subfigure}{0.3\textwidth}
    \includegraphics[width=\textwidth]{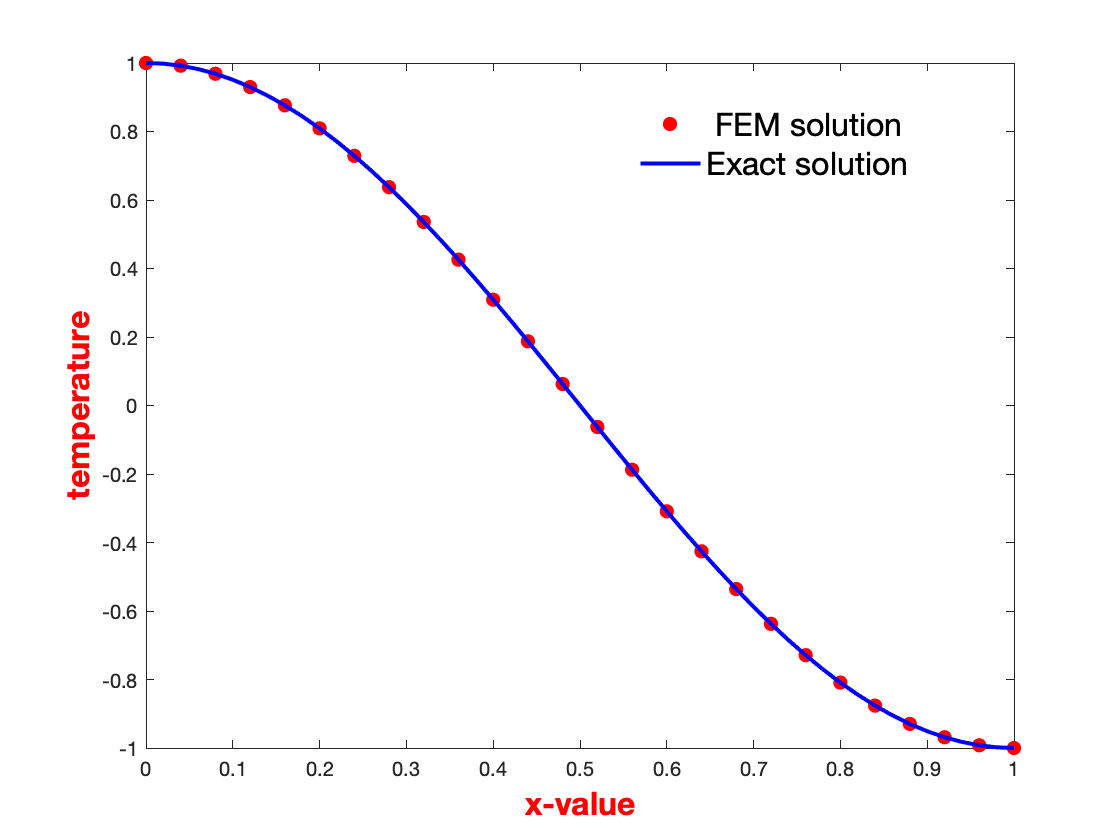}
    \caption{Plot of $\theta_h$}
    \label{fig:third}
\end{subfigure}
\caption{Comparison of the finite element solution and analytical solution for \eqref{sm}.}
\label{fig:figures}
\end{figure}

In Figure~\ref{fig:figures}, a comparison is presented for the numerical solution and the exact solution, along the midline of the domain,  as in \eqref{exact_u_p} and \eqref{exact_theta}. In the actual finite element computations, 
the mesh consists of $16384$ active cells with a total number of degrees of freedom $148739$ out of which $132098$ DOFs are for velocity and $16641$ DOFs are for pressure.   

Tables~\ref{tab_uh}, ~\ref{tab_wh}, ~\ref{tab_pv}, and  ~\ref{tab_theta_h} presents the exact solution and the numerical solution for each of the primitive variables. The computations done to obtain the data presented in these tables are for a fixed Rayleigh number,  i.e. $Ra=100$. The computational mesh is fixed with $64$-square elements in each direction. It is clear from each table that our finite element discretization algorithm using $\mathbb{Q}_2\mathbb{Q}_1$ Taylor-Hood pair of elements and the corresponding implementation yielded a highly accurate numerical solution.

\begin{table}[H]
\centering
\caption{Comparison of the exact solution and the numerical solution for $u$ along the line $z=0.5$ and $0 \leq x \leq 1$} 
\label{tab_uh}
\begin{tabular}{|c | c | c |  c |} 
 \hline
 x & Exact Solution ($u^{e} = \sin \pi \, x$) & $u_h$ & Error $= |  u^{e} - u_h | $ \\ [0.5ex] 
 \hline\hline
 0              &          0     &          0             &          0 \\
 \hline
0.1    &    0.309016994374947   &    0.309017   &  5.62505259038559e-09 \\
 \hline
0.2     &   0.587785252292473   &    0.587784   &  1.25229247316394e-06 \\
 \hline
0.3    &    0.809016994374947    &   0.809018   & 1.00562505256363e-06 \\
 \hline
 0.4   &     0.951056516295154    &   0.951056   &  5.16295153518698e-07 \\  
 \hline
 0.5         &               1      &        1           &            0 \\ 
 \hline
  0.6     &   0.951056516295154    &   0.951056  &    5.1629515362972e-07 \\ 
 \hline
   0.7   &     0.809016994374947    &   0.809018  &  1.00562505256363e-06 \\ 
 \hline
   0.8     &   0.587785252292473    &   0.587784   &  1.25229247327496e-06 \\ 
 \hline
   0.9     &    0.309016994374948    &    0.309017    &  5.62505247936329e-09\\ 
 \hline
   1    & 1.22464679914735e-16  &   1.22465e-16   &  3.20085264688233e-22 \\ 
 \hline
\end{tabular}
\end{table}

\begin{table}[H]
\centering
\caption{Comparison of the exact solution and the numerical solution for $u$ along the line $z=0.5$ and $0 \leq x \leq 1$} 
\label{tab_wh}
\begin{tabular}{|c | c | c |  c |} 
 \hline
 x & Exact Solution ($w^{e} = - \pi \, z  \cos \pi \, x$) & $w_h$ & Error $= |  w^{e} - w_h | $ \\ [0.5ex] 
 \hline\hline
 0    &      -1.5707963267949   &     -1.5708   &  3.67320510341607e-06\\
 \hline
0.1    &     -1.49391608237078    &    -1.49479   &   0.000873917629222287  \\
 \hline
0.2    &     -1.27080092307881    &   -1.27226  &    0.00145907692118508 \\
 \hline
0.3     &   -0.923290915245228    &   -0.925315   &    0.00202408475477167 \\
 \hline
 0.4    &    -0.485402759681367    &   -0.488064 &     0.00266124031863335 \\  
 \hline
 0.5  &    -9.61835346860895e-17   &   -0.00189136   &   0.0018913599999999  \\ 
 \hline
  0.6   &      0.485402759681366    &   0.482743  &     0.0026597596813665  \\ 
 \hline
   0.7    &     0.923290915245228   &    0.921265  &   0.00202591524522822  \\ 
 \hline
   0.8     &     1.27080092307881   &     1.26934  &    0.00146092307881496  \\ 
 \hline
    0.9      &     1.49391608237078   &     1.49305 &    0.000866082370777788 \\ 
 \hline
    1      &     1.5707963267949    &     1.5708    &  3.67320510341607e-06 \\ 
 \hline
\end{tabular}
\end{table}

\begin{table}[H]
\centering
\caption{Comparison of the exact solution and the numerical solution for $p$ along the line $x=0.5$ and $0 \leq z \leq 1$} 
\label{tab_pv}
\begin{tabular}{|c | c | c |  c |} 
 \hline
 z & Exact Solution ($p^{e} = \cos  \pi  z$) & $p_h$ & Error $= |  p^{e} - p_h | $ \\ [0.5ex] 
 \hline\hline
 0             &           1      &    1.0004   &  0.000399999999999956 \\
 \hline
0.1   &     0.951056516295154  &      0.951164 &   0.000107483704846478   \\
 \hline
0.2      &   0.809016994374947    &    0.809185 &    0.000168005625052592 \\
 \hline
0.3     &   0.587785252292473    &     0.58791 &   0.000124747707526907  \\
 \hline
 0.4    &    0.309016994374947    &    0.309051   &  3.40056250525689e-05  \\  
 \hline
 0.5   &  6.12323399573677e-17  &  -1.30469e-09   & 1.30469006123234e-09   \\ 
 \hline
  0.6     &  -0.309016994374947    &   -0.309051   &  3.40056250526799e-05  \\ 
 \hline
   0.7  &     -0.587785252292473   &     -0.58791  &   0.000124747707527018 \\ 
 \hline
   0.8    &   -0.809016994374947  &      -0.809185 &   0.000168005625052703   \\ 
 \hline
    0.9   &    -0.951056516295154  &     -0.951164   &  0.000107483704846478 \\ 
 \hline
    1         &              -1    &     -1.0004   &  0.000399999999999956 \\ 
 \hline
\end{tabular}
\end{table}

\begin{table}[H]
\centering
\caption{Comparison of the exact solution and the numerical solution for $\theta$ along the line $z=0.5$ and $0 \leq x \leq 1$} 
\label{tab_theta_h}
\begin{tabular}{|c | c | c |  c |} 
 \hline
 x & Exact Solution ($\theta^{e} = \sin ( \pi x +  \pi  z)$) & $\theta_h$ & Error $= |  \theta^{e} -  \theta_h | $ \\ [0.5ex] 
 \hline\hline
 0            &            1      &         1        &               0  \\
 \hline
 0.1    &    0.951056516295154    &    0.951056   &  5.1629515362972e-07  \\
 \hline
0.2   &     0.809016994374947   &     0.809018  &  1.00562505256363e-06 \\
 \hline
0.3      &  0.587785252292473    &    0.587784   &  1.25229247327496e-06   \\
 \hline
   0.4    &    0.309016994374948    &    0.309017   &  5.62505247936329e-09 \\  
 \hline
0.5   &   1.22464679914735e-16   &  -4.55889e-09   &  4.55889012246468e-09    \\ 
 \hline
0.6   &    -0.309016994374947    &   -0.309017    &  5.6250527014079e-09    \\ 
 \hline
  0.7  &     -0.587785252292473  &     -0.587784   &  1.25229247305292e-06   \\ 
 \hline
   0.8  &     -0.809016994374947    &   -0.809018   &  1.00562505267465e-06    \\ 
 \hline
  0.9    &   -0.951056516295154    &   -0.951056   &  5.16295153518698e-07    \\ 
 \hline
  1        &               -1       &       -1       &                0   \\ 
 \hline
\end{tabular}
\end{table}

\subsection{Comparison with previous studies}
In this second numerical exercise to test the validity of the discretization method developed with the code, we consider the natural convection in a unit square enclosure with heating done at left wall and the cold temperature maintained at the right wall. The top and bottom walls were kept thermally insulated.  The normal velocity is taken to be zero on all walls. We compute the \textit{Nusselt number} at the left boundary by using the formula \eqref{Nu_formula} with the exchange of $x$ and $z$.  The Nusselt number computed from our algorithm is compared with some previous studies on Darcy-B\`enard convection with heating from the left side for the Rayleigh numbers listed in the following table.

\begin{center}
\begin{tabular}{ |c|c|c|c| } 
\hline
$Ra$ & $Nu^{ref}$  & Present study \\ 
\hline
\hline
25  & 1.3682 \cite{hickox1981numerical} & 1.38523335 \\ 
\hline
50 & 1.9794 \cite{hickox1981numerical}  & 2.06974958\\ 
\hline
100 & 3.1018 \cite{hickox1981numerical,vsarler2000natural}  & 3.21446095 \\ 
\hline
1000 & 13.529 \cite{kosec2008local}  & 13.83837009 \\ 
\hline
10000  & 44.295 \cite{kosec2008local} & 43.72971012 \\ 
\hline
\end{tabular}
\end{center}

From the above table it is clear that the Nusselt number obtained in this study are ``close'' to the values obtained in \cite{kosec2008local,vsarler2000natural,hickox1981numerical}. The discrepancies in the values are believed to be due to several reasons such as the different approximation methods, diverse computational meshes, and varied linear algebra tools for solving a system of linear equations arising from a particular discretization technique. Nevertheless, the comparison presented in the above table is encouraging.   \\

The following figure \ref{fig2_DBC_left} is obtained for different Rayleigh numbers on a fixed computational mesh with $256 \times 256$ square elements. The following is the mesh information: 
\begin{center}
\begin{tabular}{ |l|l| } 
\hline
Number of active cells & 65536\\
\hline
Total number of DOFs & 592387\\ 
\hline
Velocity DOFs & 526338 \\ 
\hline
Pressure DOFs & 66049 \\ 
\hline
Temperature DOFs & 526338 \\ 
\hline
\end{tabular}
\end{center}

\begin{figure}[H]
\centering
\begin{subfigure}{0.4\textwidth}
    \includegraphics[width=\textwidth]{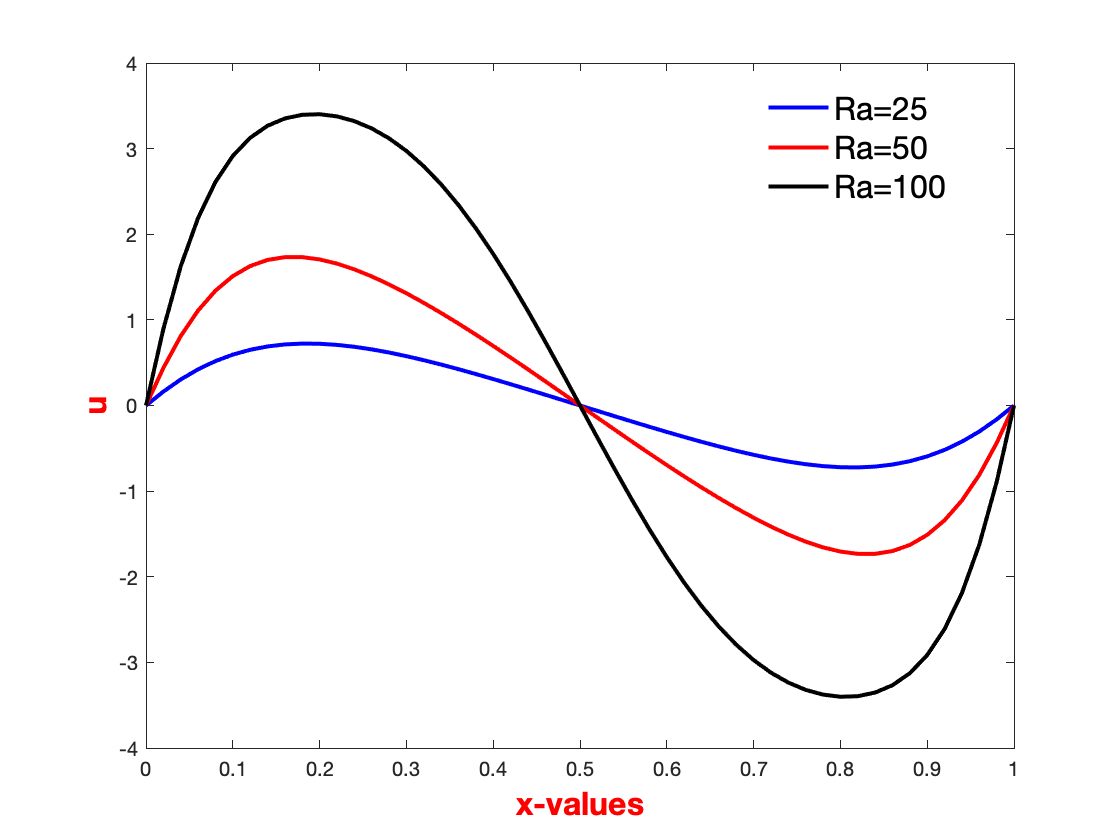}
    \caption{Plot of $u$ on the horizontal line}
\end{subfigure}
\quad 
\begin{subfigure}{0.4\textwidth}
    \includegraphics[width=\textwidth]{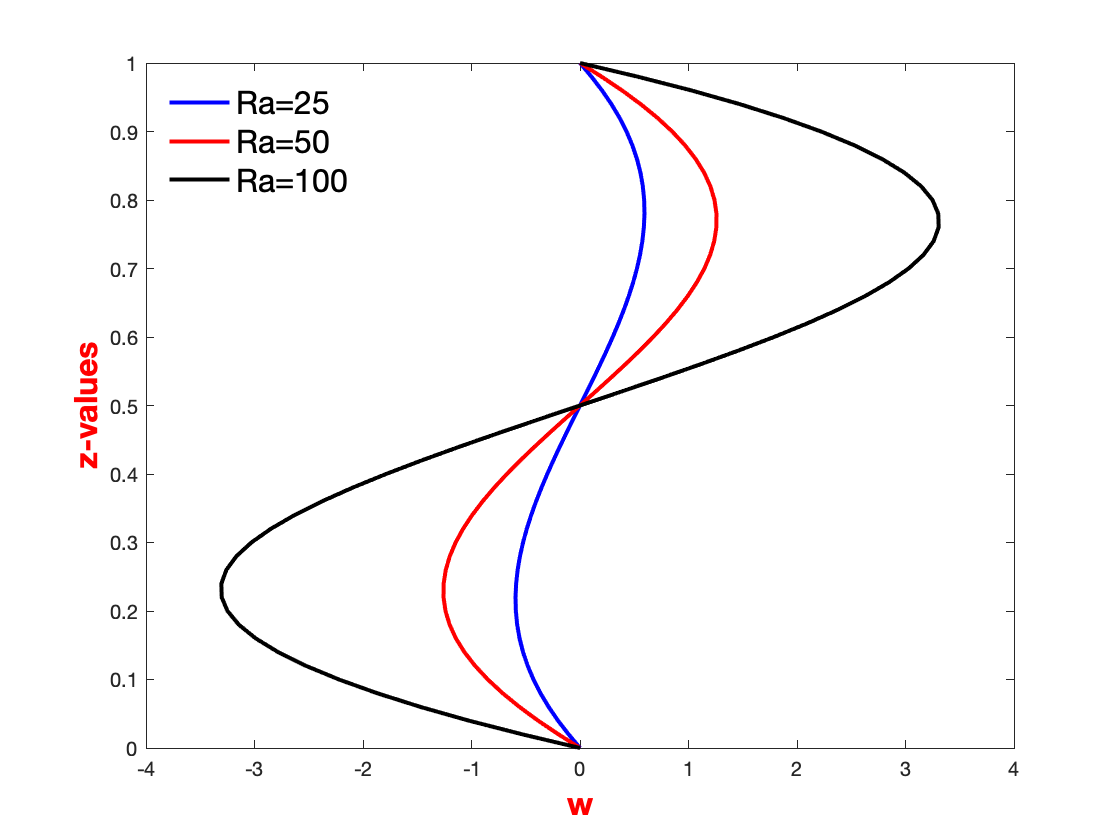}
    \caption{Plot of $w$ on the vertical line} 
\end{subfigure}
\quad 
\begin{subfigure}{0.4\textwidth}
    \includegraphics[width=\textwidth]{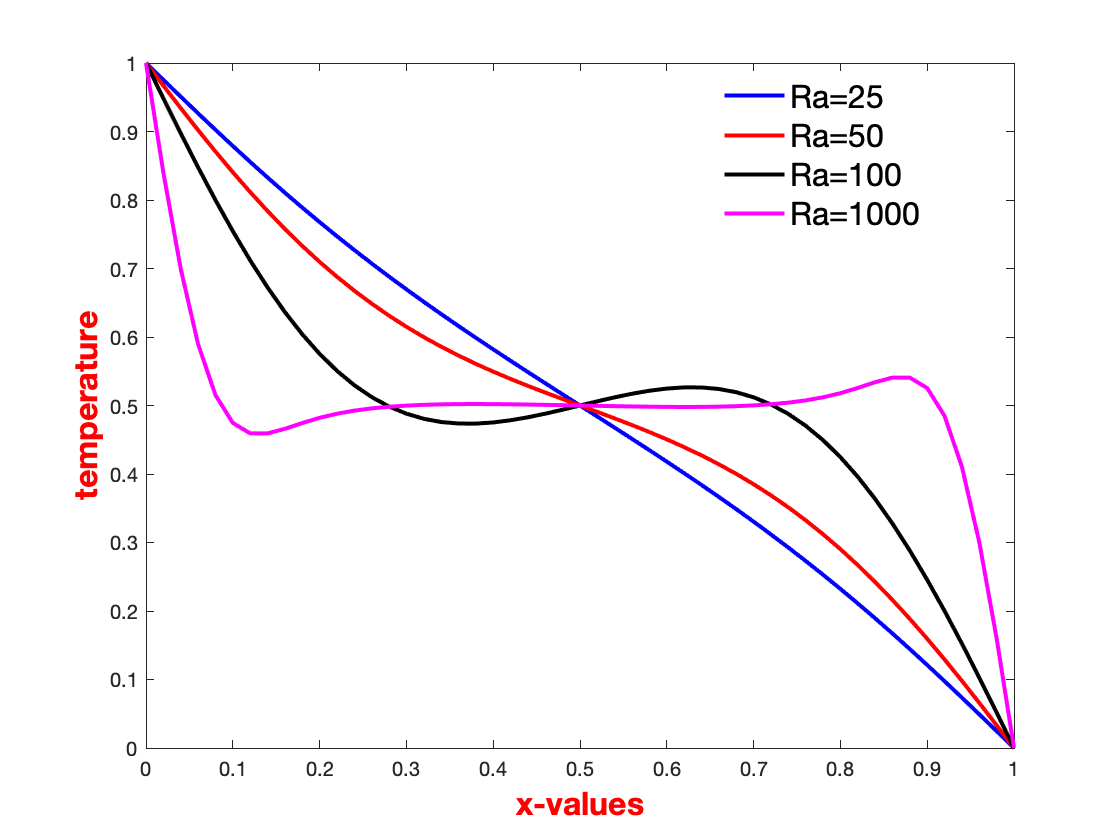}
    \caption{Plot of $\theta$ on the horizontal line}
\end{subfigure}
\quad 
\begin{subfigure}{0.4\textwidth}
    \includegraphics[width=\textwidth]{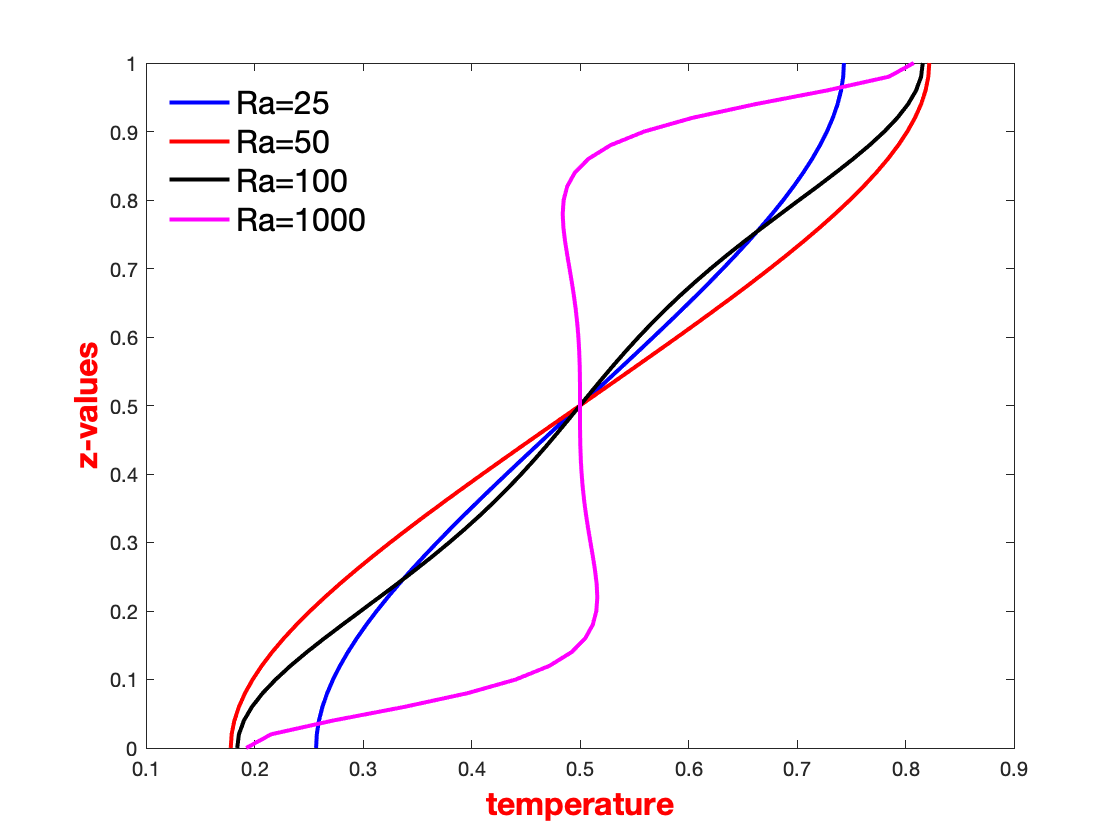}
    \caption{Plot of $\theta$ on the vertical line}
\end{subfigure}
\caption{Solution variables for the DBC problem with heating from left-wall and cooling from right-wall.  }
\label{fig2_DBC_left}
\end{figure}


\subsection{Effect of mesh density}
The main aim of this numerical test is to depict the effect of the mesh density on the primitive variables such as velocities, pressure, and temperature,  and also on the heat transfer rate, and average Nusselt number. We consider the steady-state model \eqref{sm} and its corresponding finite element formulation of the discrete problem. Using the Taylor-Hood elements, the numerical model is implemented for $5$ different meshes with the uniform mesh size of  $h=1/16, \; h=1/32, \; h=1/64, \; h=1/128, \; h=1/256$. The Rayleigh number, $Ra=100$,  is kept constant for all the $5$ numerical simulations while the hydraulic resistivity and thermal conductivity were taken as $1$. 

\begin{table}[H]
\centering
\caption{Values of the average Nusselt number for different meshes.} 
\label{tab_Nu_diff_meshes}
\begin{tabular}{ |l|l| }
\hline
$h$ & $Nu$\\
\hline
$1/16$ & 3.2572568\\ 
\hline
$1/32$ & 3.2751202 \\ 
\hline
$1/64$ & 3.29223171 \\ 
\hline
$1/128$ & 3.30579158 \\ 
\hline
$1/256$ & 3.31446095 \\ 
\hline
\end{tabular}
\end{table}

The table \ref{tab_Nu_diff_meshes} shows the value of Nusselt number for various mesh sizes. It is clear that the value of $Nu$ is somewhat converging when the mesh size is sufficiently smaller.

\begin{figure}[H]
\centering
\begin{subfigure}{0.4\textwidth}
    \includegraphics[width=\textwidth]{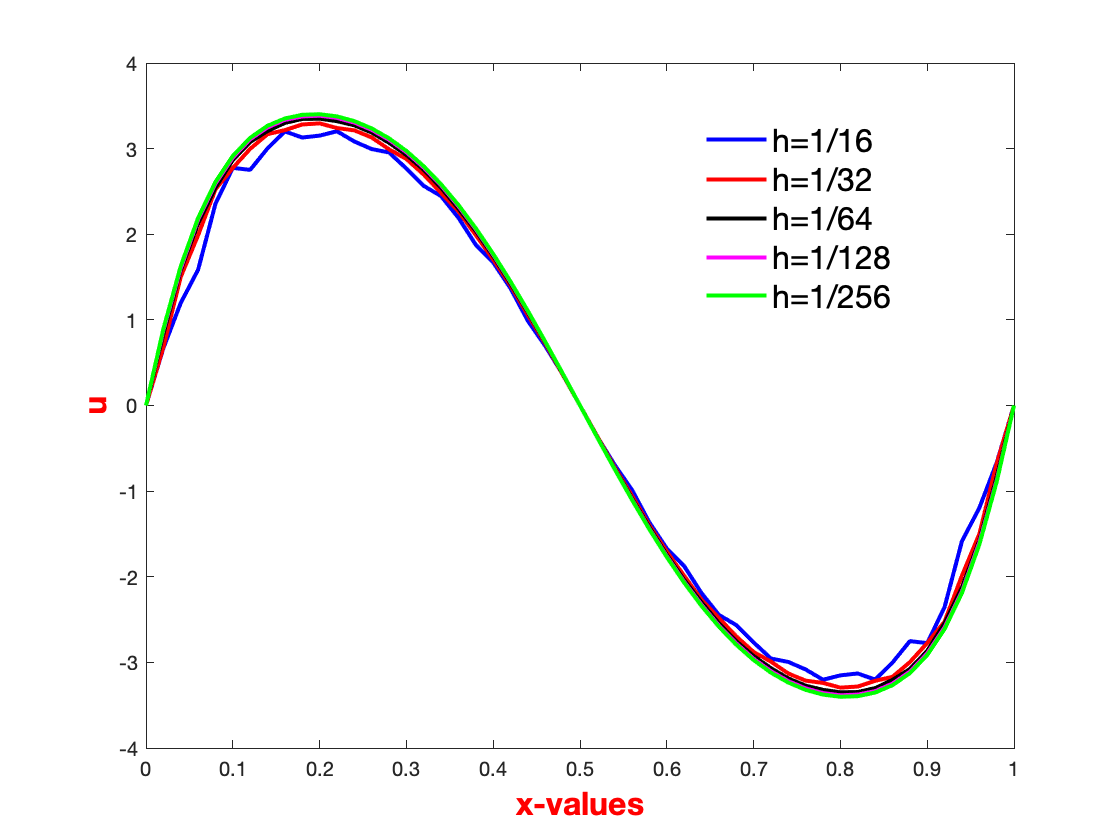}
    \caption{Plot of $u$ on the horizontal line}
\end{subfigure}
\quad 
\begin{subfigure}{0.4\textwidth}
    \includegraphics[width=\textwidth]{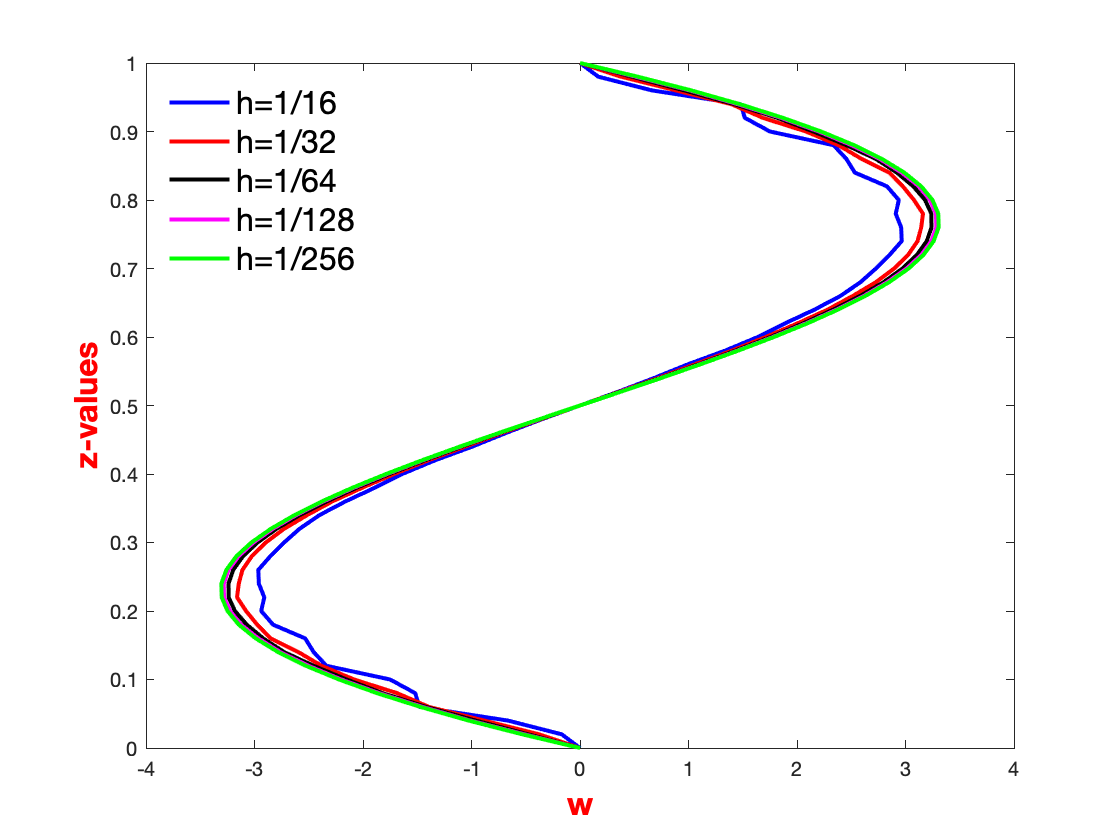}
    \caption{Plot of $w$ on the vertical line} 
\end{subfigure}
\quad 
\begin{subfigure}{0.4\textwidth}
    \includegraphics[width=\textwidth]{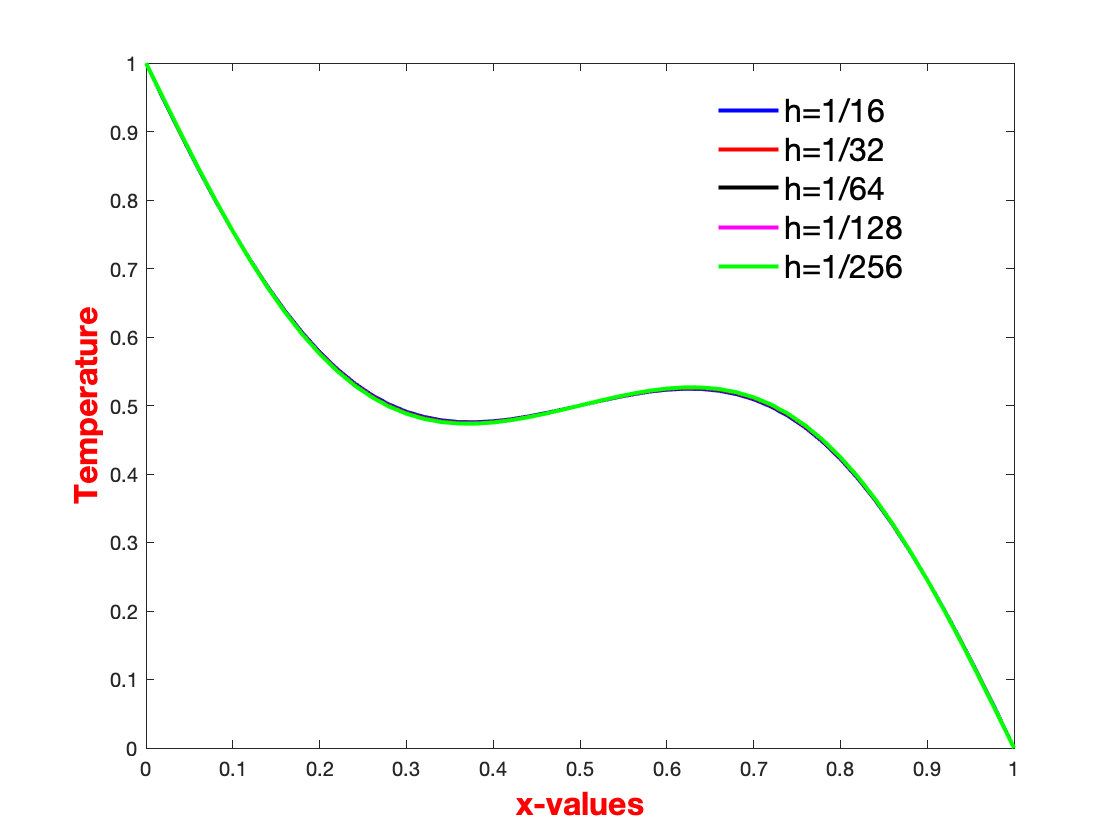}
    \caption{Plot of $\theta$ on the horizontal line}
\end{subfigure}
\quad 
\begin{subfigure}{0.4\textwidth}
    \includegraphics[width=\textwidth]{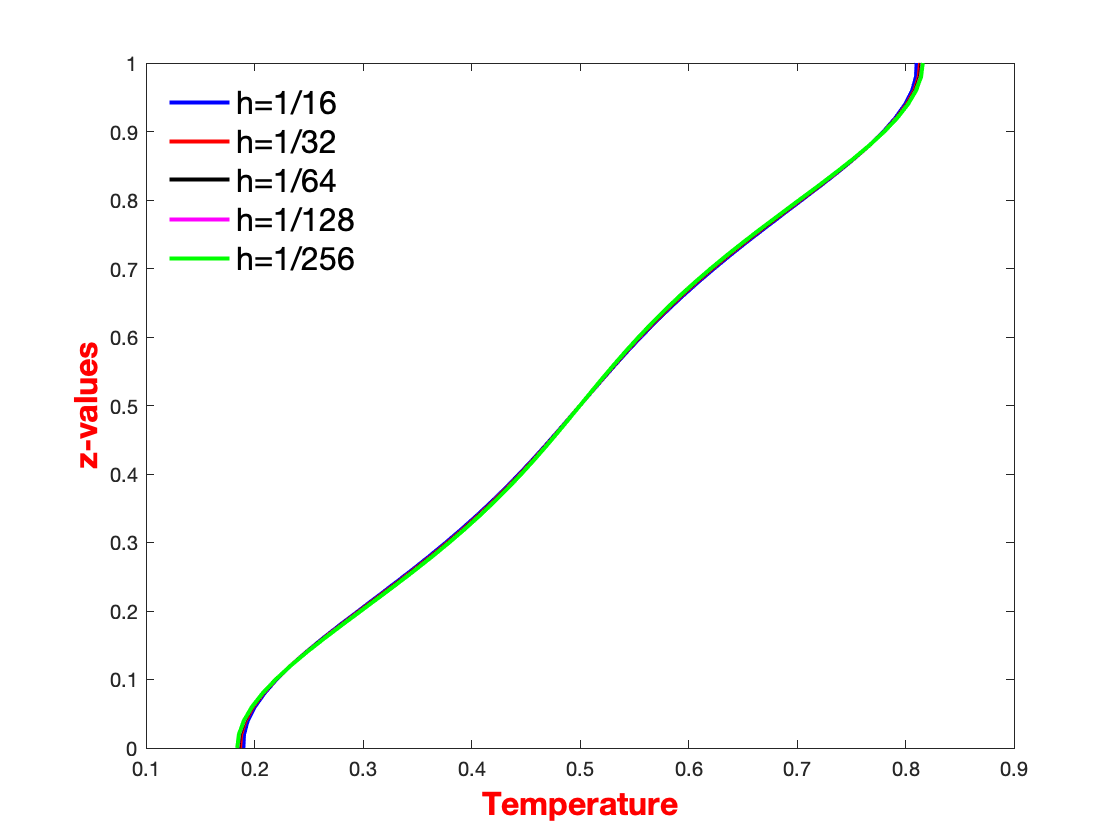}
    \caption{Plot of $\theta$ on the vertical line}
\end{subfigure}
\caption{Solution variables for the DBC problem with heating from left-wall and cooling from right-wall.  }
\label{fig4_left_diff_h}
\end{figure}

Figure~\ref{fig4_left_diff_h} depicts the solution variables, $\bfa{u}_h$ and $\theta_h$,  along the mid-lines in the domain for various mesh sizes. It is clear that some sort of convergence is seen for the smaller mesh sizes of $h=1/128$ and $h=1/256$. By looking at the numerical results, we choose $h = 1/64$ for all future computations. 

\subsection{Unsteady thermal convection}
For this numerical test, we consider the domain $\Omega = [0, \; 1] \times  [0, \; 1]$ and discretize it into $64 \times 64$ Taylor-Hood elements for velocity-pressure and bi-quadratic elements to approximate the temperature variable. Such a discretization can smoothly approximate the solution variations near the boundaries where the solution is expected to have a lot of variation. We perform the following two cases of computations based on the nature of hydraulic resistivity and thermal diffusivity:
\begin{itemize}
\item[(I)] $\chi(z) = 1+z$ and  $\zeta(z) = 1+z$,
\item[(II)] $\chi(z) = 1+z+z^2$ and  $\zeta(z) = 1+z+z^2$.
\item[(III)] $\chi(z) = 1+z$ and  $\zeta(z) = 1+z+z^2$,
\item[(IV)] $\chi(z) = 1+z+z^2$ and  $\zeta(z) = 1+z$.
\end{itemize}
The above choices for the functions $\chi$ and $\zeta$ are made to illustrate some test cases, however, our computational code is made to accept any arbitrary choices for these two functions. The linear and quadratic variations are preferred in porous media convection problems \cite{bhatta2017convective,bhatta2020thermal}. 
In each case, the numerical calculations are done for various  Rayleigh numbers, such as $Ra = 50, \; 100, \; 500, \; 1000$,  with a uniform bottom boundary heating while the rest of the boundaries are kept sealed. No-slip condition is assumed for the velocity variable on all the boundaries.  For each of the above cases, our goal is to characterize the effect of $Ra$ and provide the plots for velocity contours along with isotherm lines. The numerical results have also been obtained for the local Nusselt numbers for various Rayleigh numbers and a detailed study is presented for each case considered above. 

\subsubsection{Case-I: Linear variations in  both $\chi(z)$ and $\zeta(z)$ and the effect of $Ra$}
\begin{figure}[H]
\centering
\begin{subfigure}{0.4\textwidth}
    \includegraphics[width=\textwidth]{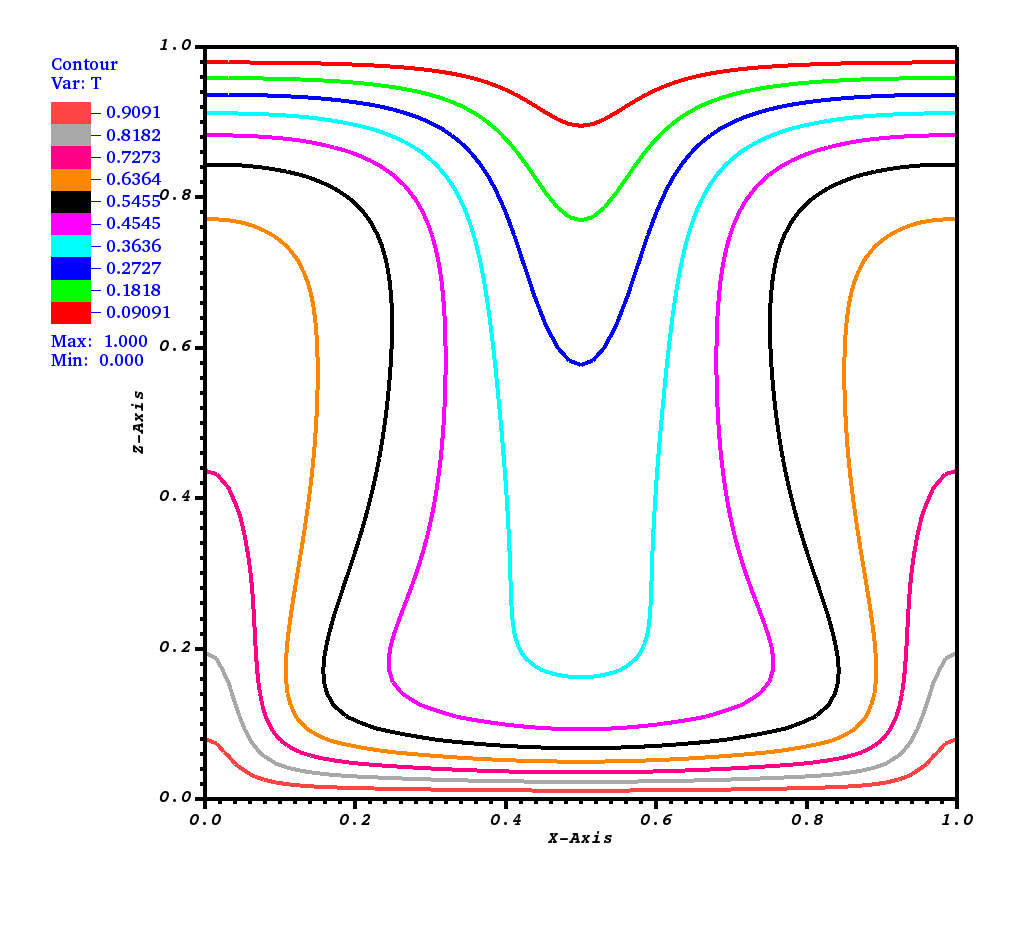}
    \caption{Isotherms for $Ra=500$}
\end{subfigure}
\quad 
\begin{subfigure}{0.4\textwidth}
    \includegraphics[width=\textwidth]{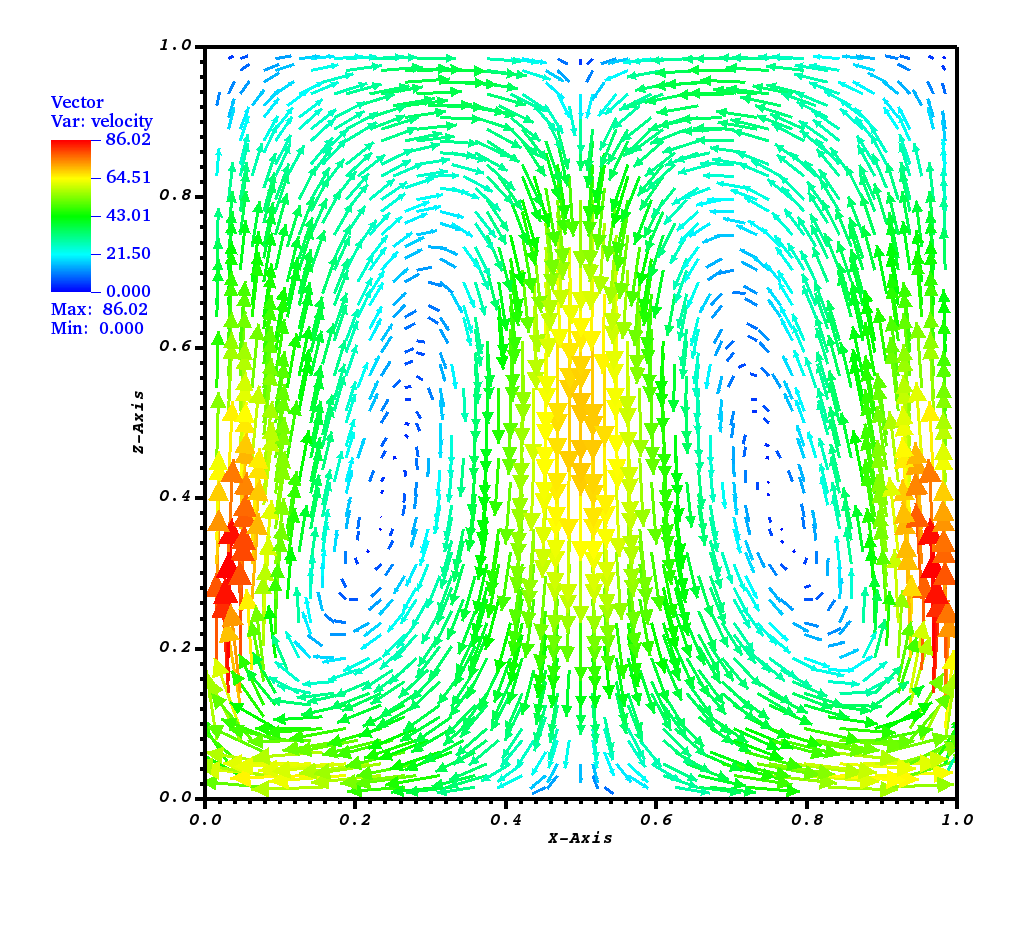}
    \caption{$\bfa{u}$-vector for $Ra=500$} 
\end{subfigure}
\quad 
\begin{subfigure}{0.4\textwidth}
    \includegraphics[width=\textwidth]{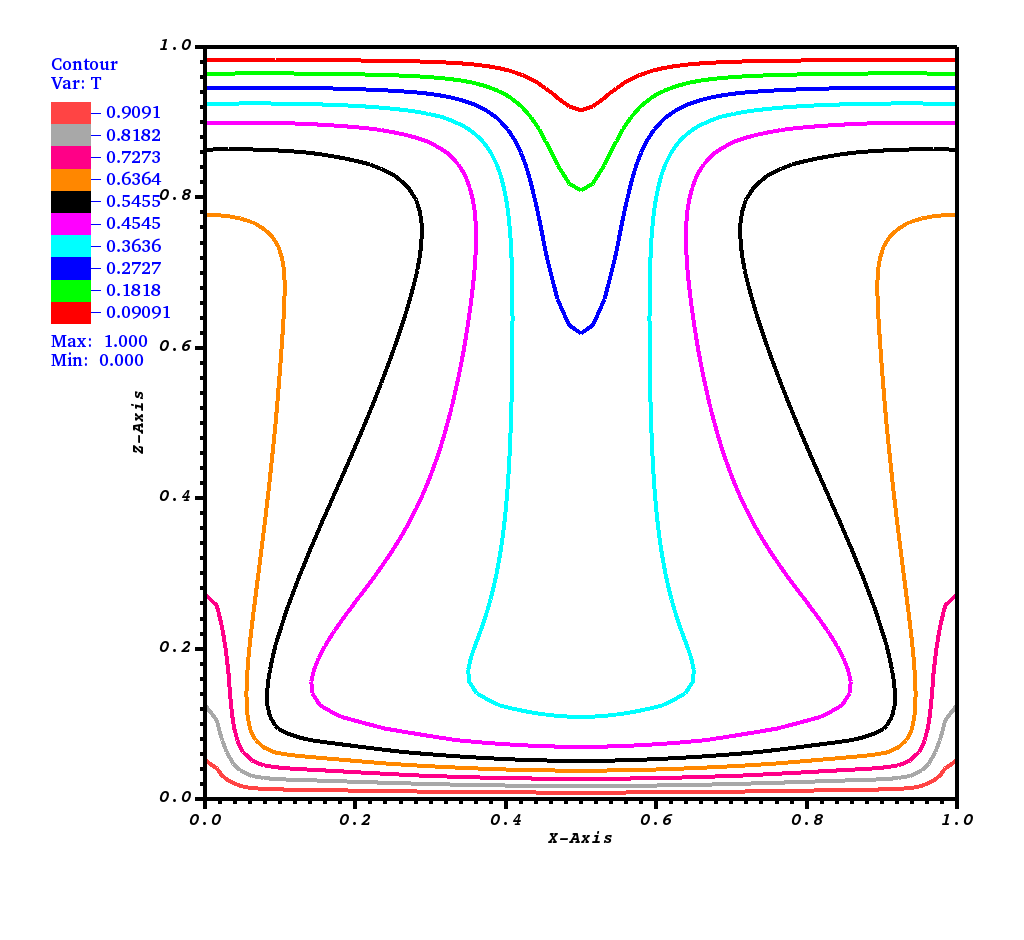}
    \caption{Isotherms for $Ra=1000$}
\end{subfigure}
\quad 
\begin{subfigure}{0.4\textwidth}
    \includegraphics[width=\textwidth]{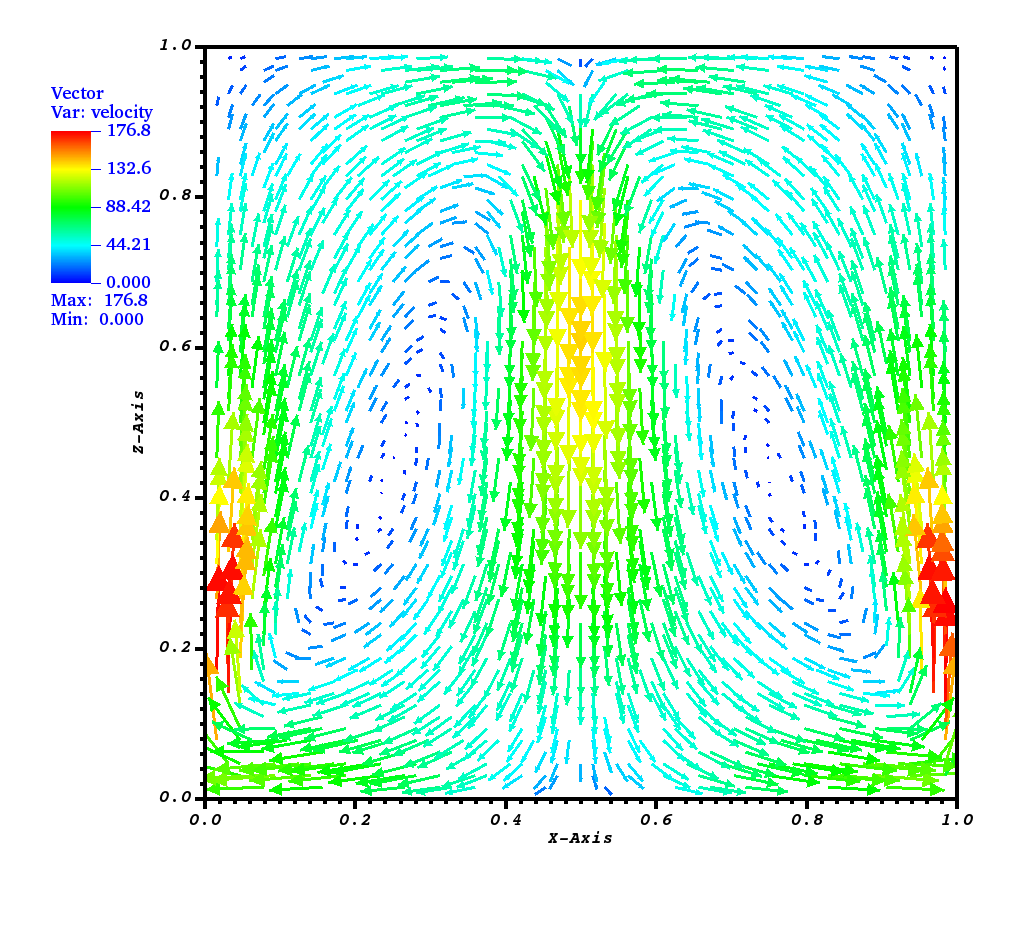}
    \caption{$\bfa{u}$-vector for $Ra=1000$}
\end{subfigure}
\caption{Results for Case-I }
\label{fig_case_I}
\end{figure}

In this test, we consider both hydraulic resistivity and thermal diffusivity as linear functions depend only on the vertical variation. This subsection also demonstrates the influence of the Rayleigh number on the flow behavior in the porous medium as well as the local heat transfer at the bottom hot plate. We analyze the numerical results based on the isotherms and magnitude of the velocity vector. The influence of $Ra$ on the flow structure is presented in Figure~\ref{fig_case_I}. The effect of $Ra$ is shown for $Ra=500, \; 1000$. The increase in $Ra$ enhances the circulations, hence creating two strong cells. Higher values of $Ra$ galvanize the circulation guided predominantly by the buoyancy, which is shown in Figure~\ref{fig_case_I}. A clear upward motion of fluid is seen near the boundaries and downward flow is prevalent at the center of the domain. 

Another interesting point is that for a fixed $Ra$, the temperature is higher at the right and left boundary as one moves from left to the center (or from the right to the center) of the domain. The temperature is always lower at the top and higher at the bottom as it should be due to the applied Dirichlet boundary conditions. As the Rayleigh number increases,  the temperature is distributed towards the wall and more predominant effects are seen. Also, for the higher $Ra$ values, the top half of the temperature is higher at the center than the lower half.  For the uniform heating of the bottom wall, Table~\ref{table_nu_case_I} presents the average values of the Nusselt numbers at the lower boundary. 
\begin{table}[H]
\centering
\caption{Values of the average Nusselt number for different meshes.} 
\label{table_nu_case_I}
\begin{tabular}{ |l|l| }
\hline
$Ra$ & $Nu$\\
\hline
$50$ & 1.4426368 \\ 
\hline
$100$ & 1.4426432 \\ 
\hline
$500$ & 6.21777759 \\ 
\hline
$1000$ & 8.55699109 \\ 
\hline
\end{tabular}
\end{table}
It is clear from the above table that the Nusselt numbers increase with increasing Rayleigh numbers. For the type of uniform heating considered in this investigation, The average Nusselt numbers remains constant up to $Ra=100$.   
\begin{figure}[H]
\centering
\includegraphics[width=0.6\textwidth]{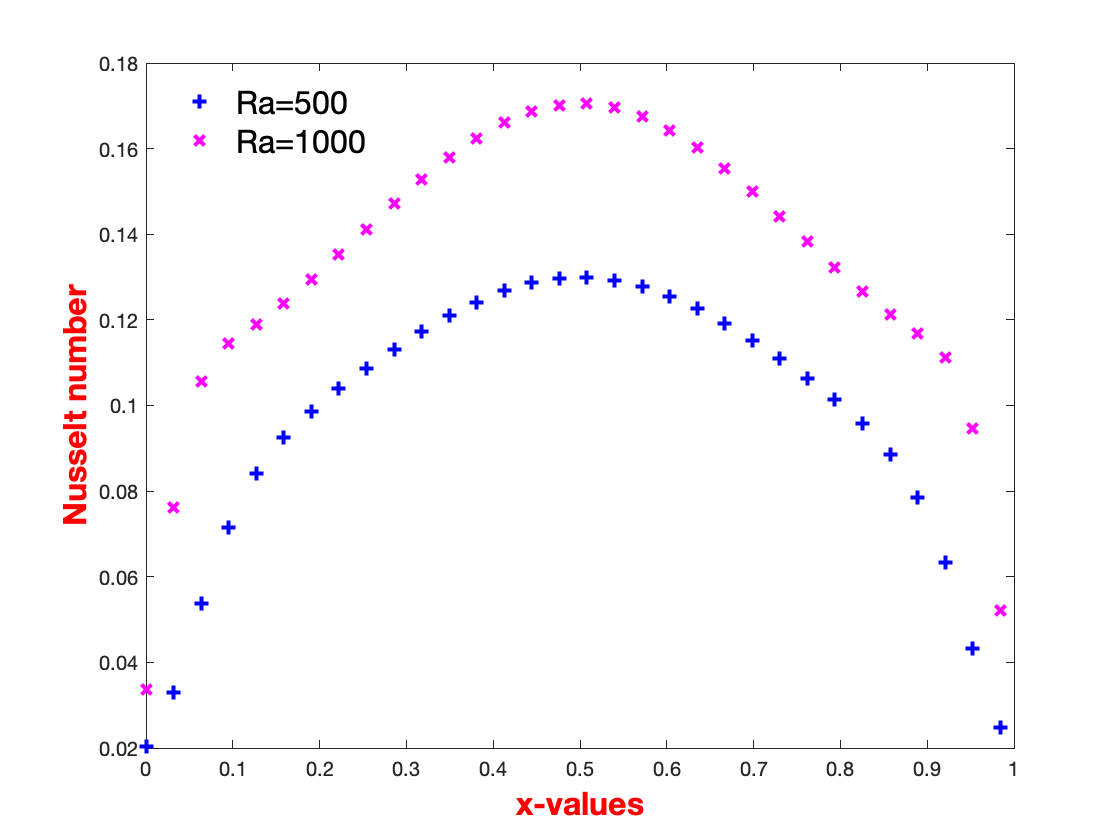}
\caption{Local Nusselt number for different values of Rayleigh number. }
\label{fig_nu_case1}
\end{figure}
In Figure~\ref{fig_nu_case1}, the Nusselt number on the lower boundary is presented for different Rayleigh numbers. For $Ra>100$, the maximum occurs at the middle of $z=0$. it may be noted from the above figure that the heat transfer rate increases from the left to the middle of the boundary and then later decreases. Such a symmetric nature of the Nusselt numbers indicates that the fluid motion will have two symmetric cells in the domain. 

\subsubsection{Case-II: Quadratic variation in both $\chi(z)$ and  $\zeta(z)$ and the effect of $Ra$}
The numerical test presented in this subsection demonstrates the influence of hydraulic resistivity and thermal diffusivity on the flow field. For the numerical experiment, we choose $\chi(z) =  \zeta(z) = 1 + z + z^2$. We also perform flow visualization for different values of Rayleigh number. 

\begin{figure}[H]
\centering
\begin{subfigure}{0.4\textwidth}
    \includegraphics[width=\textwidth]{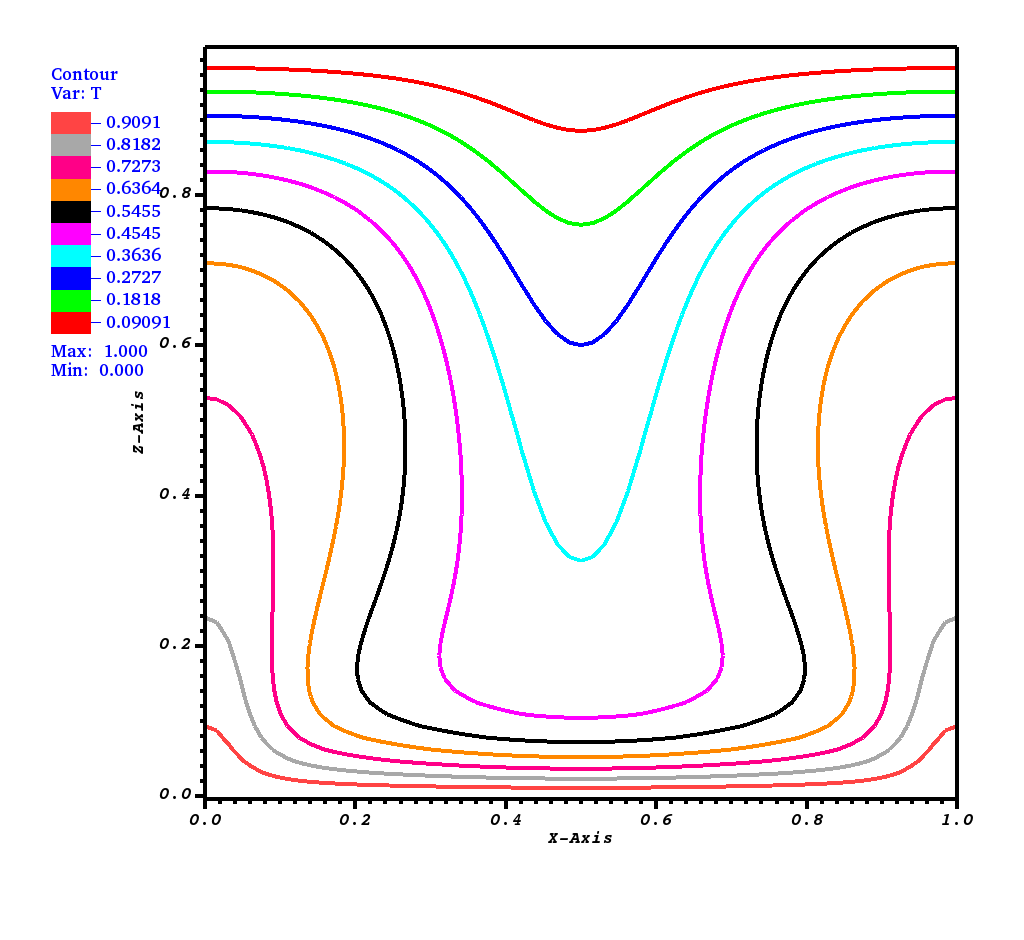}
    \caption{Isotherms for $Ra=500$}
\end{subfigure}
\quad 
\begin{subfigure}{0.4\textwidth}
    \includegraphics[width=\textwidth]{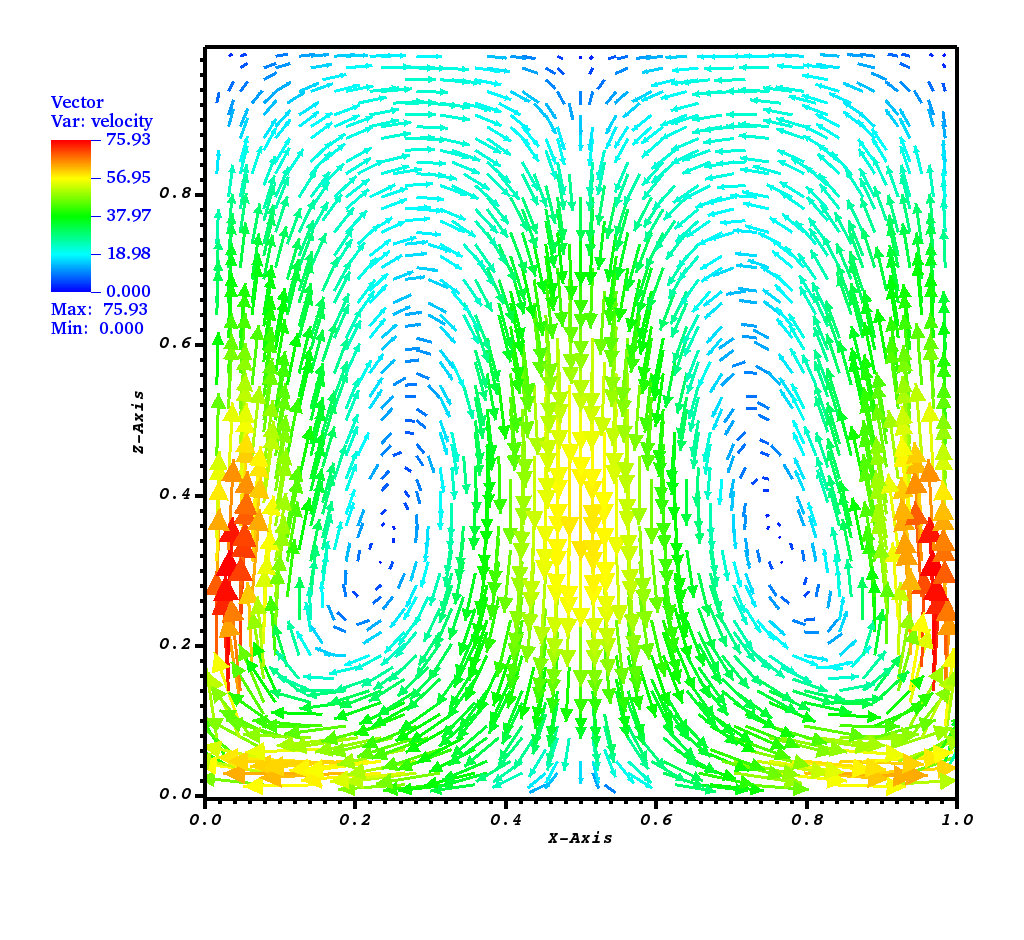}
    \caption{$\bfa{u}$-vector for $Ra=500$} 
\end{subfigure}
\quad 
\begin{subfigure}{0.4\textwidth}
    \includegraphics[width=\textwidth]{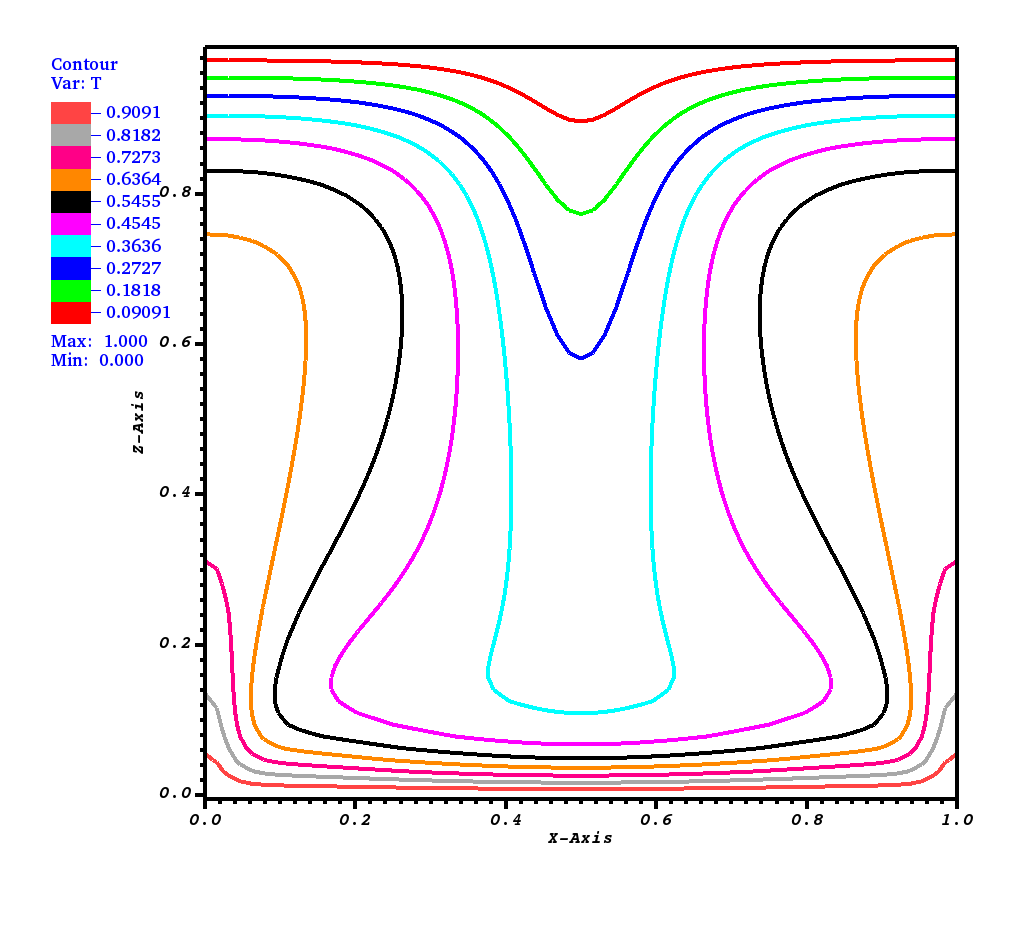}
    \caption{Isotherms for $Ra=1000$}
\end{subfigure}
\quad 
\begin{subfigure}{0.4\textwidth}
    \includegraphics[width=\textwidth]{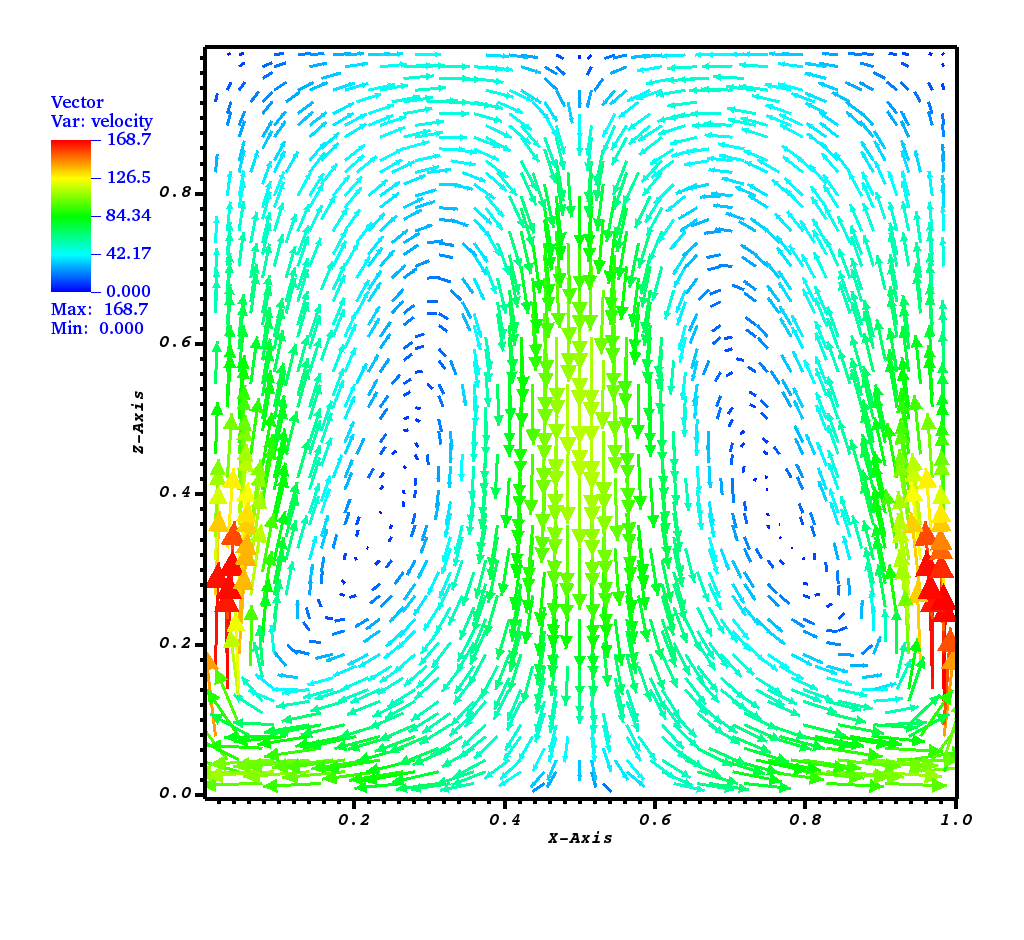}
    \caption{$\bfa{u}$-vector for $Ra=1000$}
\end{subfigure}
\caption{Results for Case-II }
\label{fig_case_II}
\end{figure}

Figure~\ref{fig_case_II} presents the effect of Rayleigh number on the flow field as well as on the temperature isotherm distribution. For the higher values of the Rayleigh number, the isotherm profile exhibit more distortion which leads to stronger and upward fluid velocity profile. Both of these two leads of higher energy transport which was clearly identified in the increase in the value of Nusselt number. 

\begin{table}[H]
\centering
\begin{tabular}{ |l|l| }
\hline
$Ra$ & $Nu$\\
\hline
$50$ & 1.65398399 \\ 
\hline
$100$ & 1.653990399 \\ 
\hline
$500$ & 5.77938259 \\ 
\hline
$1000$ & 8.6355727 \\ 
\hline
\end{tabular}
\caption{Average Nusselt number for different values of Rayleigh numbers.} 
\label{table_nu_case_II}
\end{table}

\begin{figure}[H]
\centering
\includegraphics[width=0.6\textwidth]{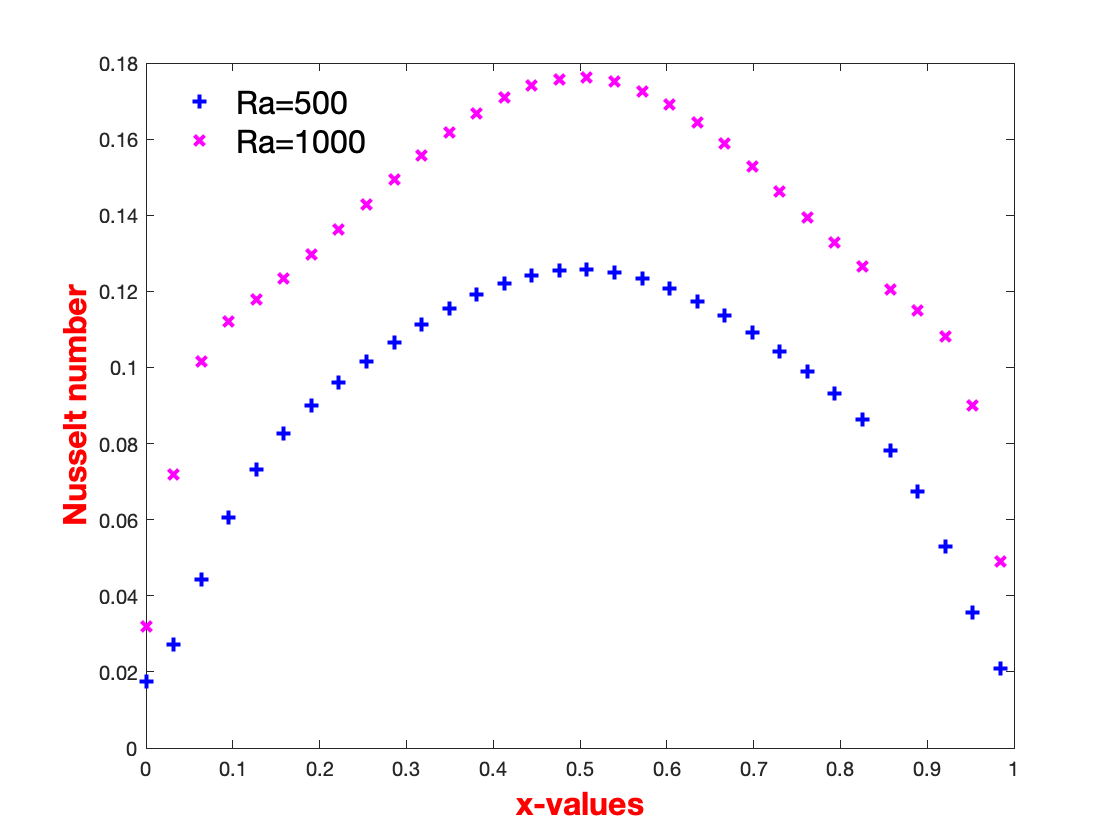}
\caption{Local Nusselt number for different values of Rayleigh number. }
\label{fig_nu_caseII}
\end{figure}

Table~\ref{table_nu_case_II} presents the values of the average Nusselt number at the bottom plate for all the Rayleigh numbers considered in this case. It is clear that the heat transfer rate increases with increasing the value of the Rayleigh number. In figure~\ref{fig_nu_caseII}, the local Nusselt number on the bottom plate for different values of Rayleigh numbers is presented. The largest value of the local Nusselt number is in the middle of the bottom edge and there is a clear symmetry in the values which indicates that there will be two cells. 

\subsubsection{Case-III: Linear variation in $\chi(z)$ and quadratic variation in  $\zeta(z)$ and the effect of $Ra$}
\begin{figure}[H]
\centering 
\begin{subfigure}{0.4\textwidth}
    \includegraphics[width=\textwidth]{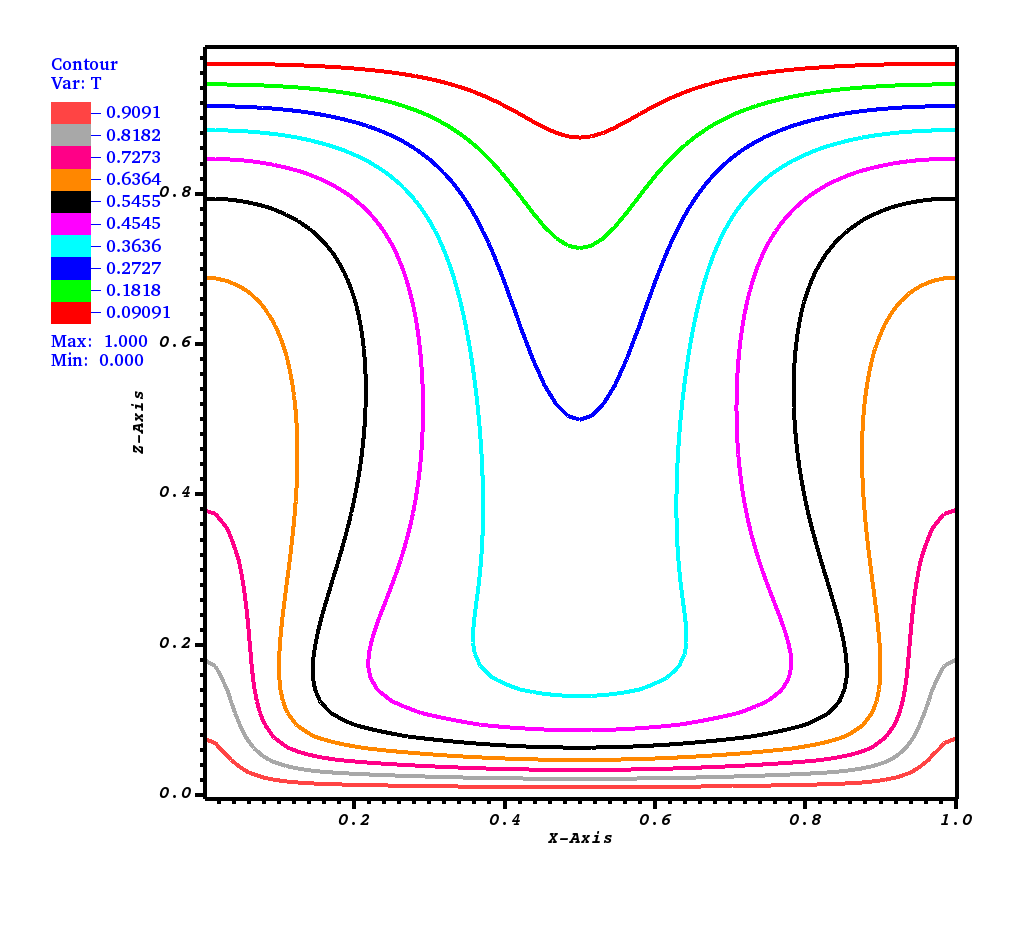}
    \caption{Isotherms for $Ra=500$}
\end{subfigure}
\quad 
\begin{subfigure}{0.4\textwidth}
    \includegraphics[width=\textwidth]{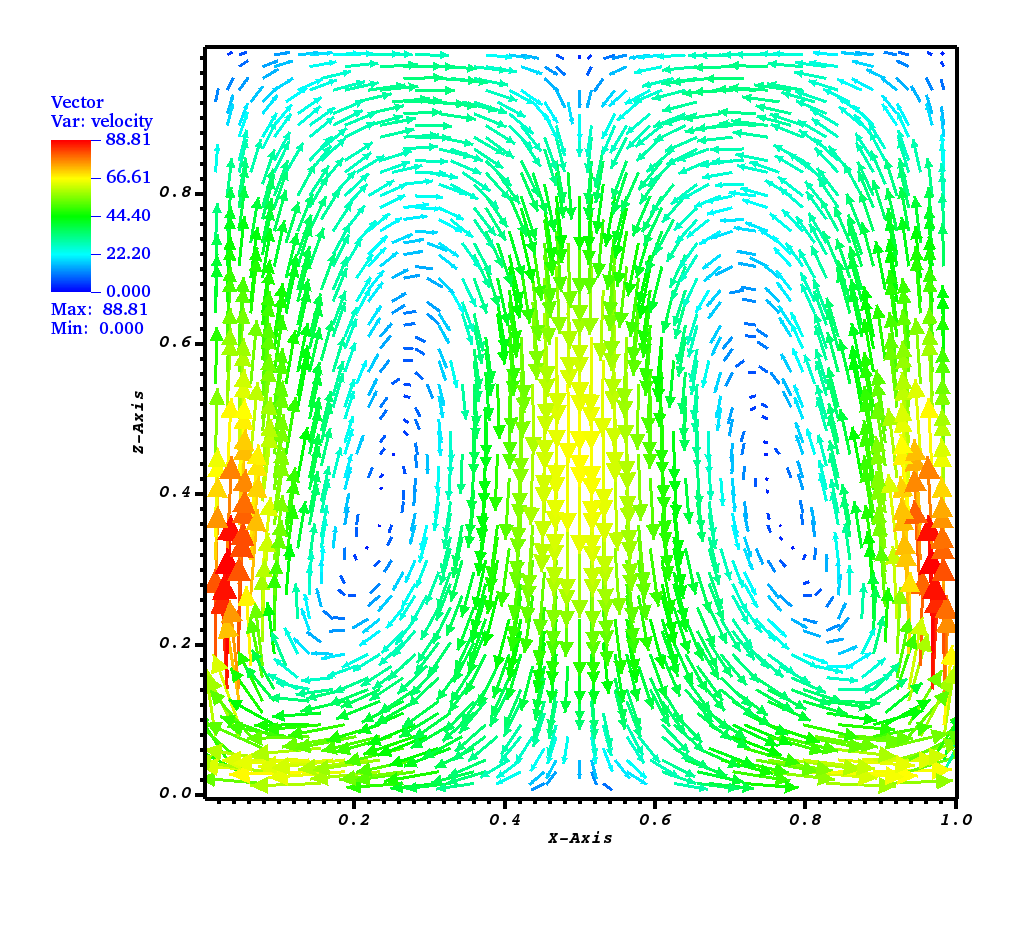}
    \caption{$\bfa{u}$-vector for $Ra=500$} 
\end{subfigure}
\quad 
\begin{subfigure}{0.4\textwidth}
    \includegraphics[width=\textwidth]{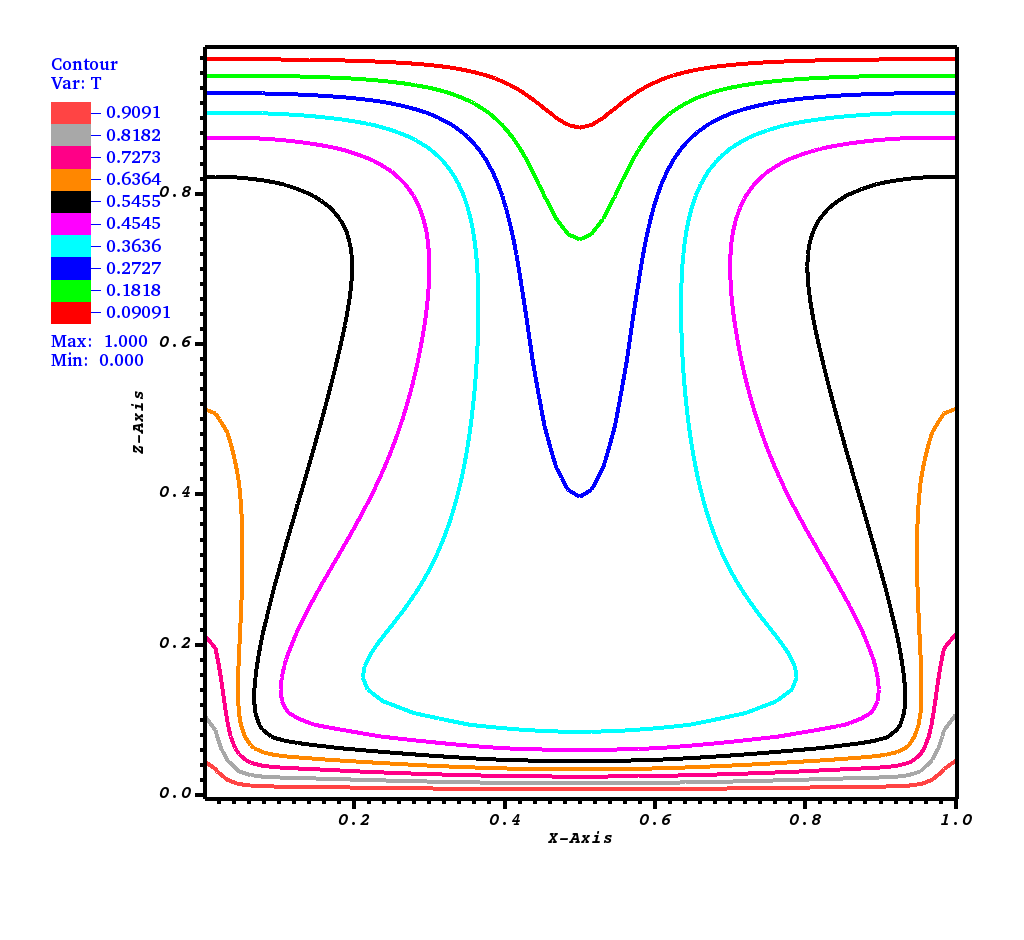}
    \caption{Isotherms for $Ra=1000$}
\end{subfigure}
\quad 
\begin{subfigure}{0.4\textwidth}
    \includegraphics[width=\textwidth]{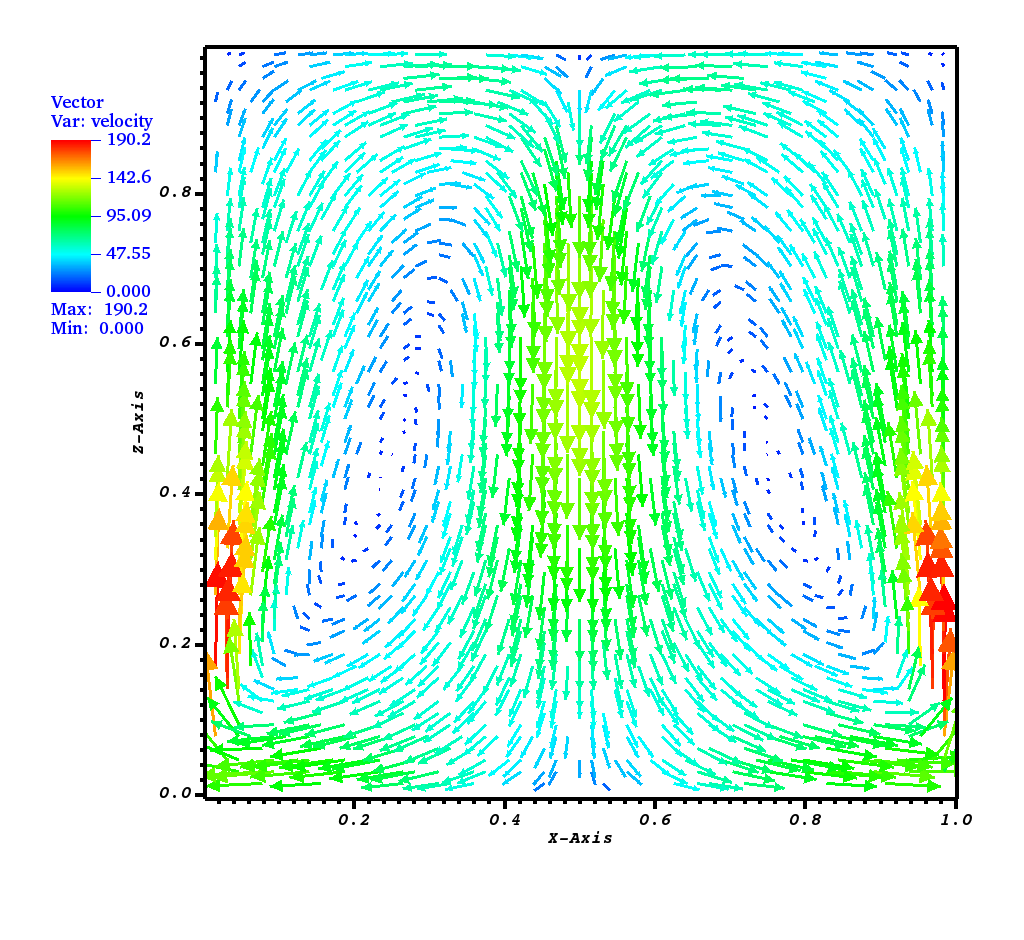}
    \caption{$\bfa{u}$-vector for $Ra=1000$}
\end{subfigure}
\caption{Results for Case-III }
\label{fig_caseIII}
\end{figure}

In this section, we consider the following types of variations in hydraulic resistivity and thermal diffusivity terms:
\begin{equation*}
\chi(z) = 1+ z, \quad \zeta(z) = 1 + z + z^2. 
\end{equation*}
Our goal is to study the effect of the Rayleigh number on thermal convection in a porous medium with the above variations in the properties. Figure~\ref{fig_caseIII} depicts both isotherms and velocity vectors for two different values of the Rayleigh number. As expected, due to quadratic variation in the thermal diffusivity term, the isotherms vary more smoothly compared to the previous two cases, which resulted in the increase in the velocity vector near the vertical walls.  The temperature distribution is higher at the side walls compared to the middle portion. Hence, the hotter fluid rises up faster near the vertical walls compared to the other parts of the domain.

Table~\ref{table_nu_case_III} presents the average Nusselt number for different values of the Rayleigh numbers. From the table, it is clear that the heat transfer rate increases with the variations in hydraulic receptivity and thermal diffusivity terms considered in this section.

\begin{table}[H]
\centering
\begin{tabular}{ |l|l| }
\hline
$Ra$ & $Nu$\\
\hline
$50$ & 1.65398399 \\ 
\hline
$100$ & 1.65399039 \\ 
\hline
$500$ & 6.544582799 \\ 
\hline
$1000$ & 9.372756199 \\ 
\hline
\end{tabular}
\caption{Average Nusselt number for different values of Rayleigh numbers.} 
\label{table_nu_case_III}
\end{table}

\begin{figure}[H]
\centering
\includegraphics[width=0.6\textwidth]{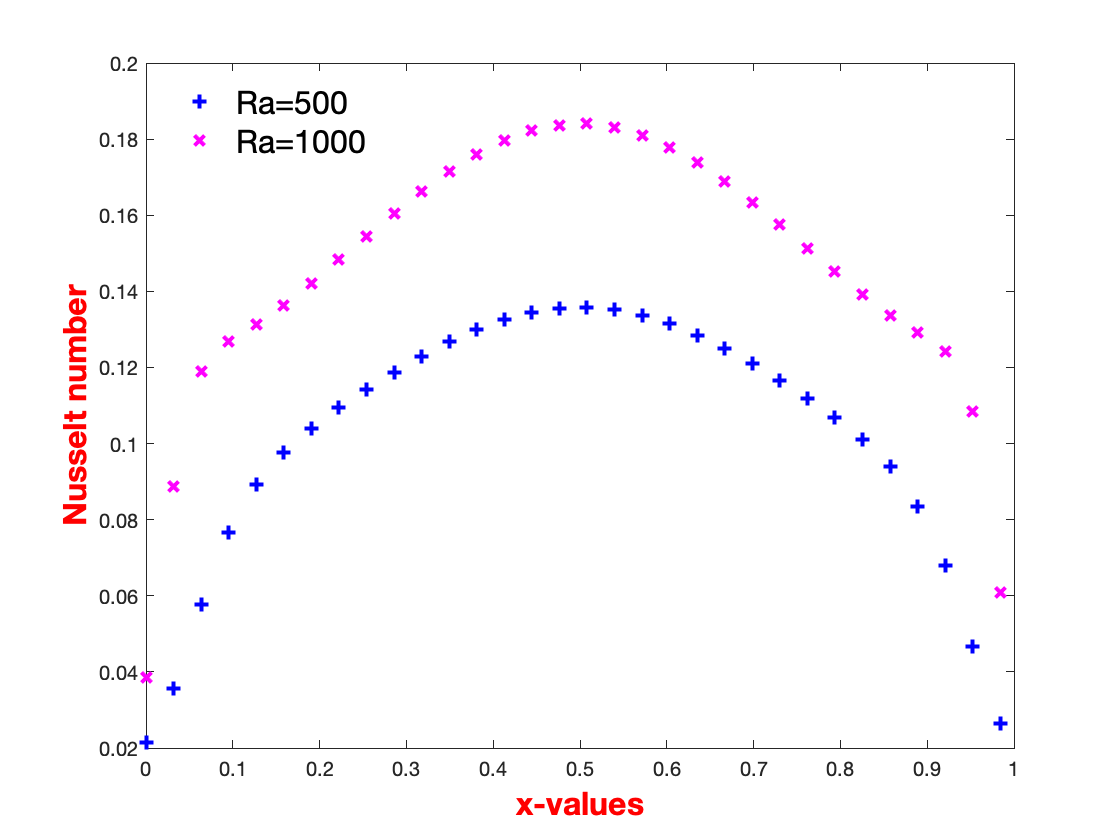}
\caption{Local Nusselt number for different values of Rayleigh number. }
\label{fig_nu_caseIII}
\end{figure}

Figure~\ref{fig_nu_caseIII} displays the local Nusselt number for different values of the Rayleigh number. In the present case, we have considered the uniform heating at the bottom wall, and for such a physical situation the largest $Nu_{loc}$ occurs at the middle of the bottom boundary. The values are symmetric, and hence there will be two fluid cells in the domain. The aforementioned trend is consistent with the previous cases studied in this paper.

\subsubsection{Case-IV: Quadratic variation in $\chi(z)$ and linear variation in  $\zeta(z)$ and the effect of $Ra$}
\begin{figure}[H]
\centering
\begin{subfigure}{0.4\textwidth}
    \includegraphics[width=\textwidth]{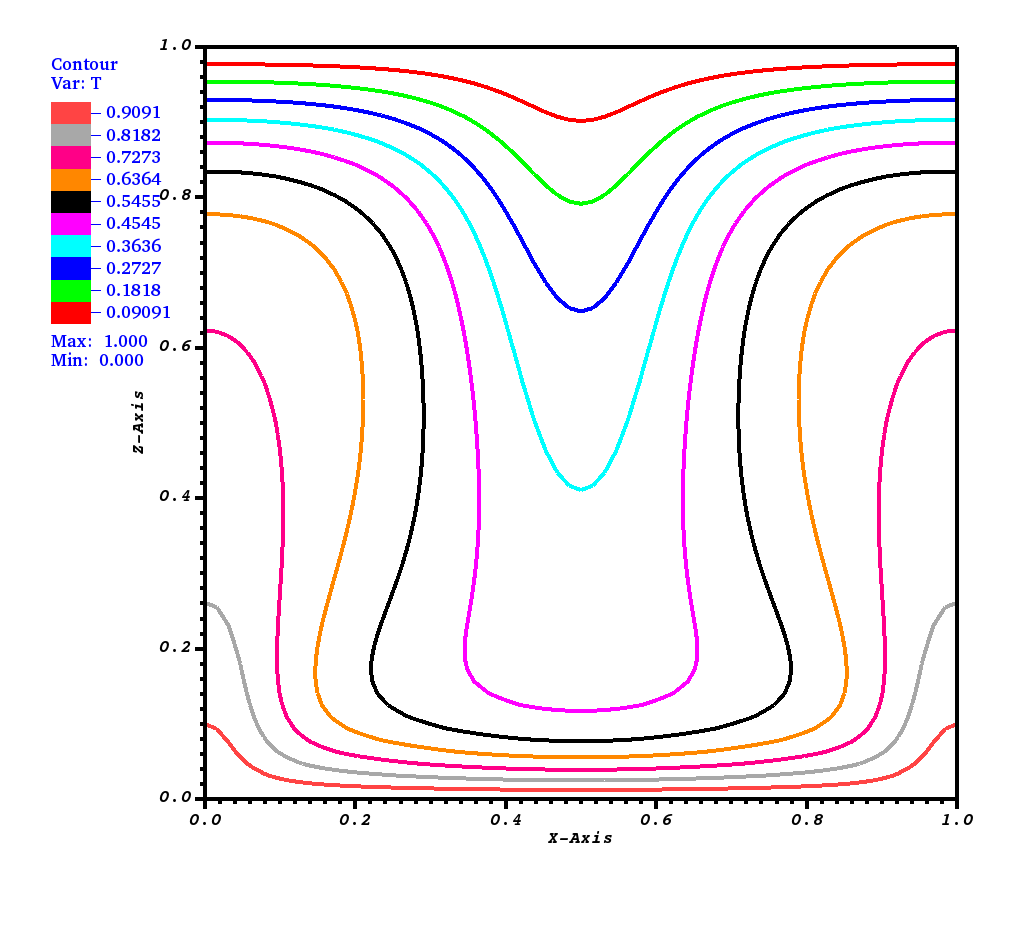}
    \caption{Isotherms for $Ra=500$}
\end{subfigure}
\quad 
\begin{subfigure}{0.4\textwidth}
    \includegraphics[width=\textwidth]{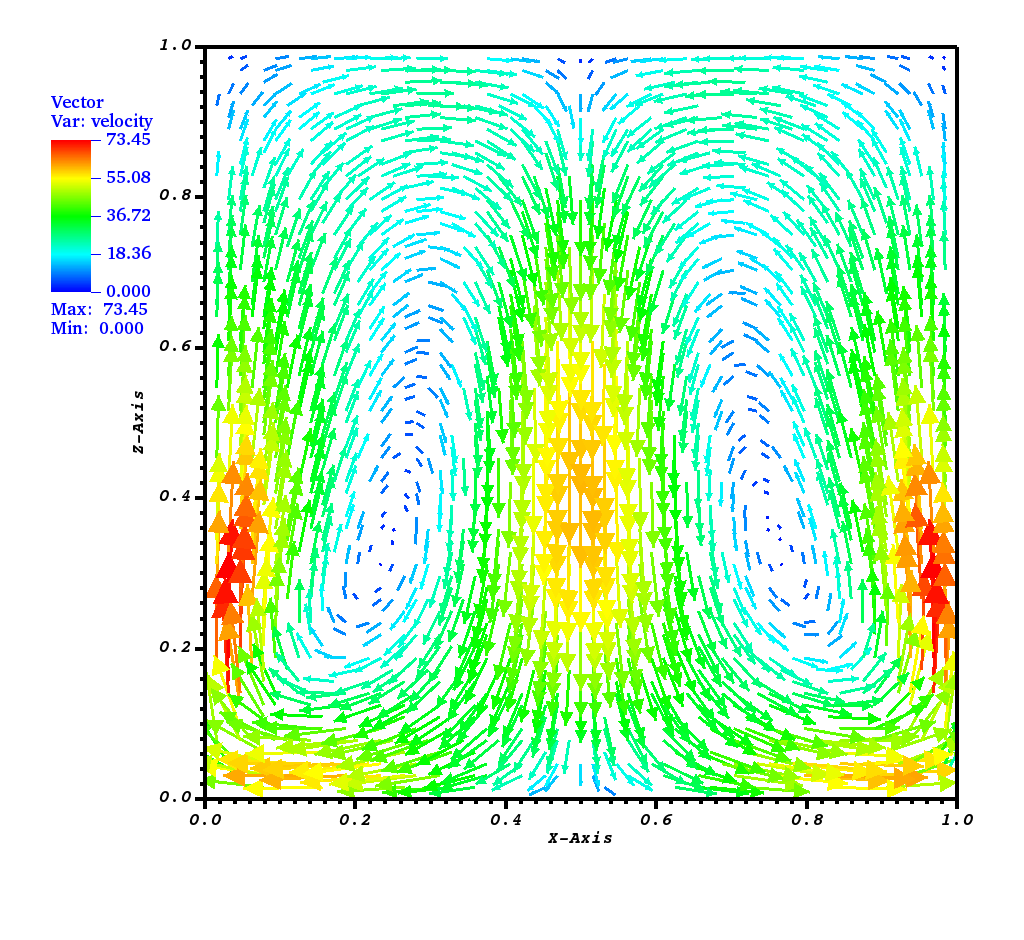}
    \caption{$\bfa{u}$-vector for $Ra=500$} 
\end{subfigure}
\quad 
\begin{subfigure}{0.4\textwidth}
    \includegraphics[width=\textwidth]{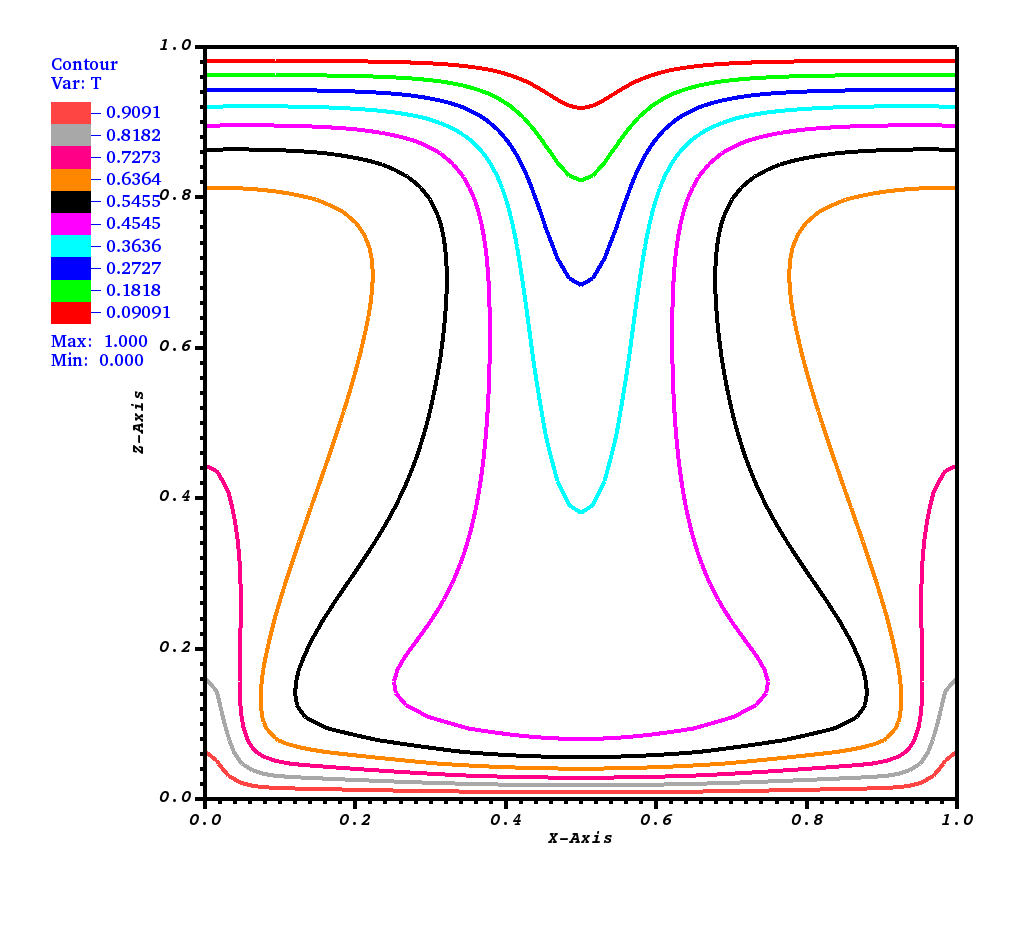}
    \caption{Isotherms for $Ra=1000$}
\end{subfigure}
\quad 
\begin{subfigure}{0.4\textwidth}
    \includegraphics[width=\textwidth]{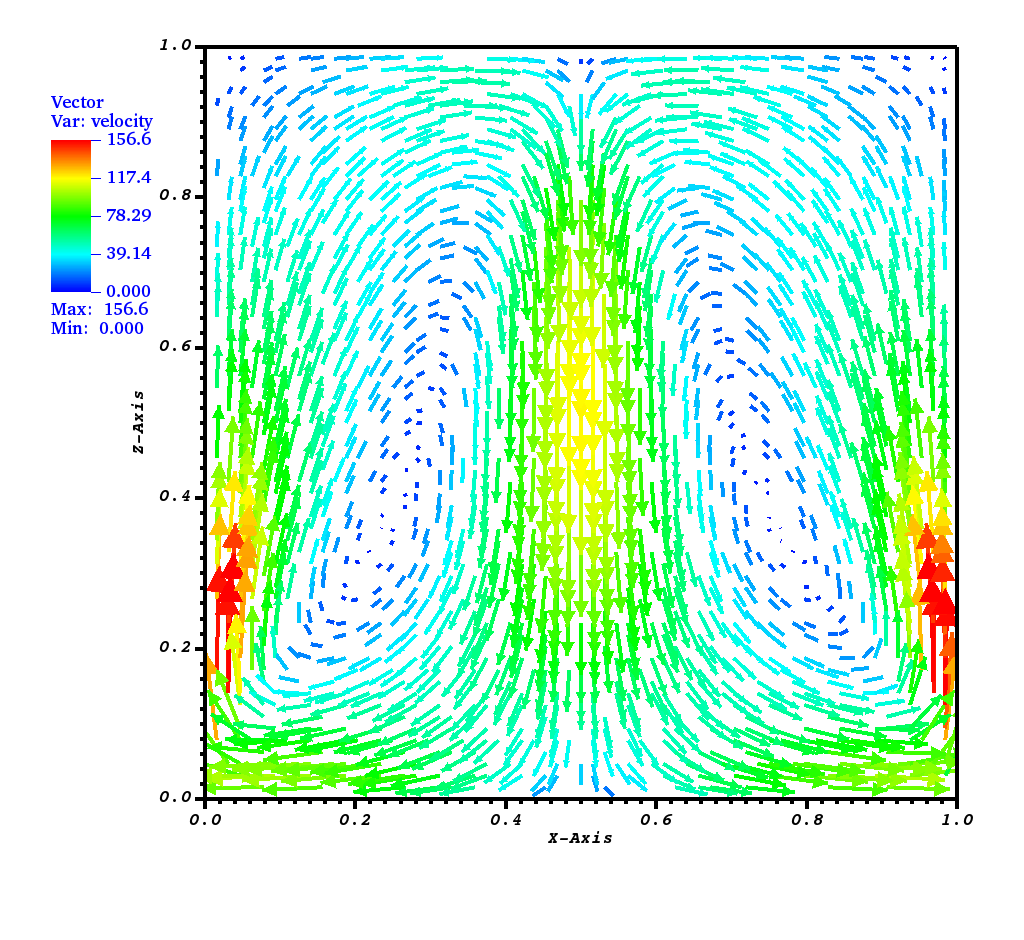}
    \caption{$\bfa{u}$-vector for $Ra=1000$}
\end{subfigure}
\caption{Results for Case-IV }
\label{fig_caseIV}
\end{figure}

In the current subsection, we consider the quadratic variation for the hydraulic resistivity term $\chi(z)$ and linear variation for the thermal diffusivity term $\zeta(z)$.  Figure~\ref{fig_caseIV} illustrates isotherms and velocity vectors for different values of the Rayleigh numbers. From the figure, it is clear that for the higher value of the Rayleigh number the isotherms for $\theta_h >0.5$ occupy a substantial region near the hot boundary wall. Hence, there will large concentration of fluid near the vertical walls, and the same can be seen in the figure.

Table~\ref{table_nu_case_IV} depicts the values of the average Nusselt numbers for different values of the Rayleigh number. There is an apparent increase in the Nusselt number values for increasing values of the Rayleigh number. An increase in $Ra$ obviously enhances the heat transfer rate and such a result is consistent with the previously reported studies. 

\begin{table}[H]
\centering
\begin{tabular}{ |l|l| }
\hline
$Ra$ & $Nu$\\
\hline
$50$ & 1.4426368 \\ 
\hline
$100$ & 1.4426432 \\ 
\hline
$500$ & 5.4844262 \\ 
\hline
$1000$ & 7.85507629 \\ 
\hline
\end{tabular}
\caption{Average Nusselt number for different values of Rayleigh numbers.} 
\label{table_nu_case_IV}
\end{table}

\begin{figure}[H]
\centering
\includegraphics[width=0.6\textwidth]{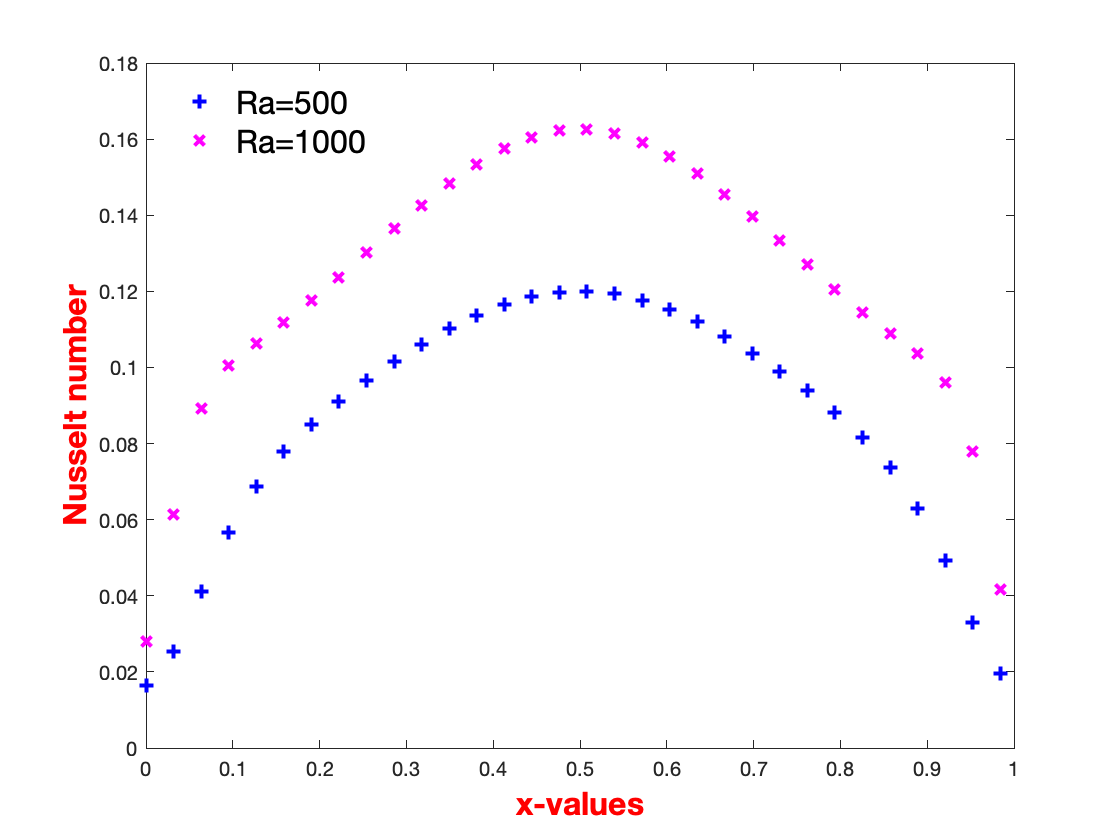}
\caption{Local Nusselt number for different values of Rayleigh number. }
\label{fig_nu_caseIV}
\end{figure}

Figure~\ref{fig_nu_caseIV} depicts the local Nusselt number on the bottom boundary of the domain. For the given domain $\Omega$ and for a given Rayleigh number, the local Nusselt number has a maximum value at the center of the bottom boundary and a symmetric behavior.

In the next part of the paper we present a detailed comparison among the four cases studied and present a conclusion based on the heat transfer rate.

\subsubsection{Comparison among different cases of variation in $\chi(z)$ and $\zeta(z)$}
In this section, we present a detailed comparison of the results obtained for all four cases of hydraulic resistivity and thermal diffusivity. For the comparison, we fixed the Rayleigh number as $Ra=500$. The effects of all the cases on the heat transfer rates along the bottom hot boundary are depicted in the table~\ref{table_nu_all_case}. It is clear that Case III, with 
\begin{equation}
\chi(z) = 1 + z, \quad \zeta(z) = 1 + z + z^2,
\end{equation}
has the highest Nusselt number. For a fixed $Ra$, the linear hydraulic resistivity and quadratic thermal diffusivity provides a better heat transfer rate compared to all other scenarios. 

\begin{table}[H]
\centering
\begin{tabular}{ | l | l | l | l | }
\hline
Cases & $\chi(z)$ & $\zeta(z)$ &  $Nu$\\
\hline
Case-I  & $1 + z$   & $1 + z$             & 6.217777599 \\ 
\hline
Case-II  & $1 + z + z^2$   & $1 + z + z^2$             & 5.77938259  \\ 
\hline
Case-III  & $1 + z$   & $1 + z + z^2$             & 6.54458279 \\ 
\hline
Case-IV  & $1 + z + z^2$   & $1 + z$             &  5.4844262 \\ 
\hline
\end{tabular}
\caption{Average Nusselt number for different values of Rayleigh numbers.} 
\label{table_nu_all_case}
\end{table}

\begin{figure}[H]
\centering
\includegraphics[width=0.6\textwidth]{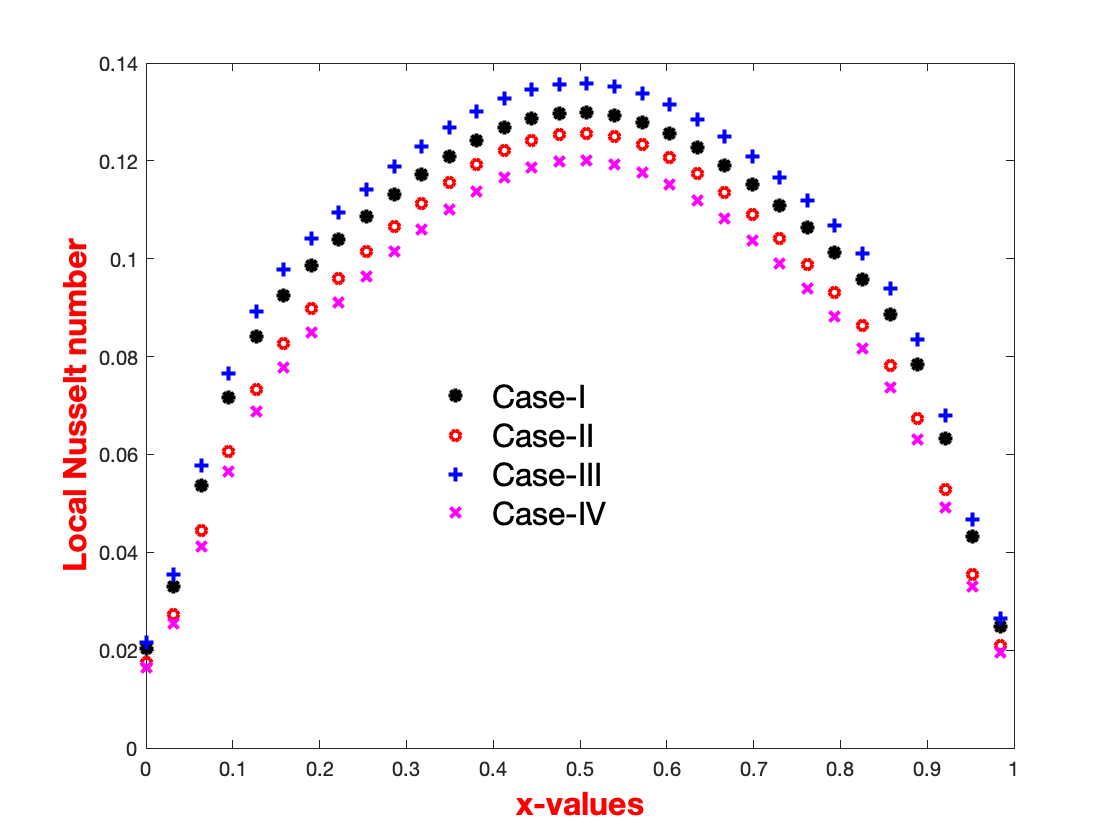}
\caption{Local Nusselt number for different variations in $\chi(z)$ and $\zeta(z)$. }
\label{fig_nu_all_cases}
\end{figure}

The overall effects of the porous medium and the fluid characteristics on the heat transfer rate at the bottom boundary are shown in figure~\ref{fig_nu_all_cases}. For case-III, displayed as the top-most curve in the figure, the Nusselt number is maximum. It is interesting to observe that the smallest value of the Nusselt number is found for the case with quadratic variation of the hydraulic resistivity and linear variation in the thermal diffusivity. 

The local Nusselt number shown in the above figure~\ref{fig_nu_all_cases} indicates that the maximum, in all the cases, occurs at the middle of the bottom boundary. The curves are symmetric about the line 
\[
x=0.5, \quad 0 \leq z \leq 1,
\]
which indicates that there will be two convection cells. 
\begin{figure}[H]
\centering
\begin{subfigure}{0.4\textwidth}
    \includegraphics[width=\textwidth]{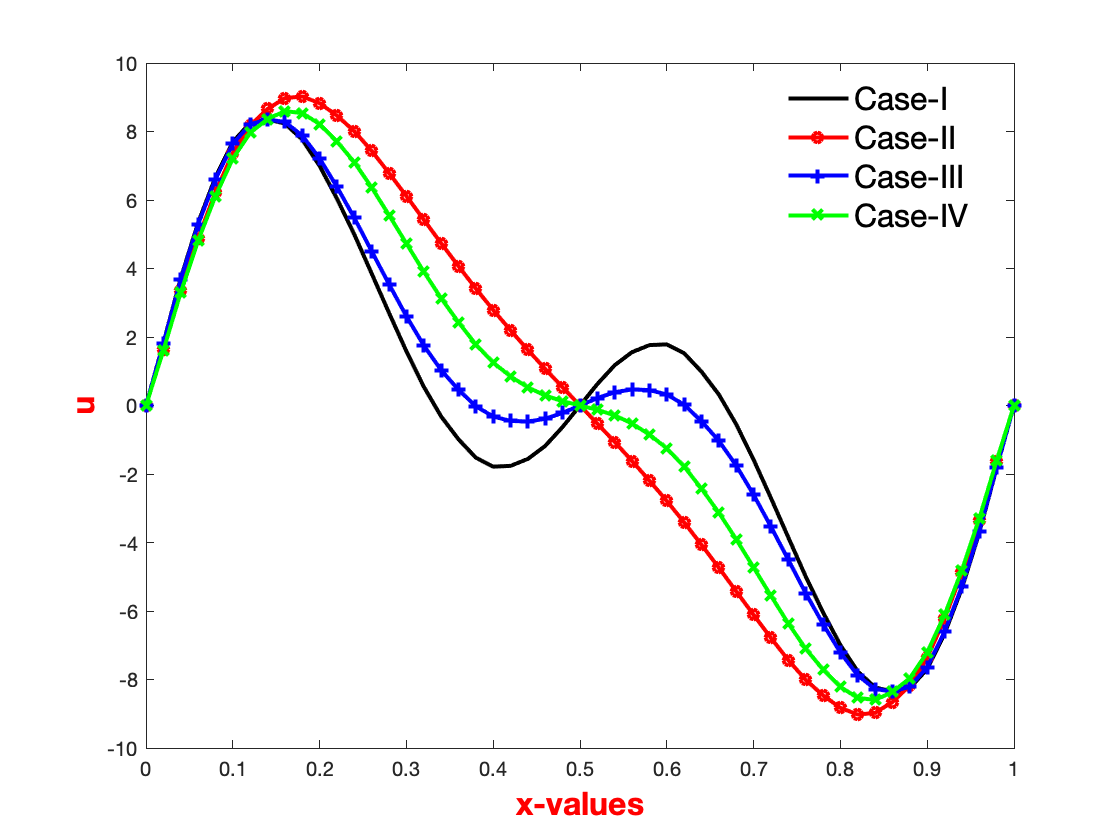}
    \caption{$u$ along horizontal midline}
\end{subfigure}
\quad 
\begin{subfigure}{0.4\textwidth}
    \includegraphics[width=\textwidth]{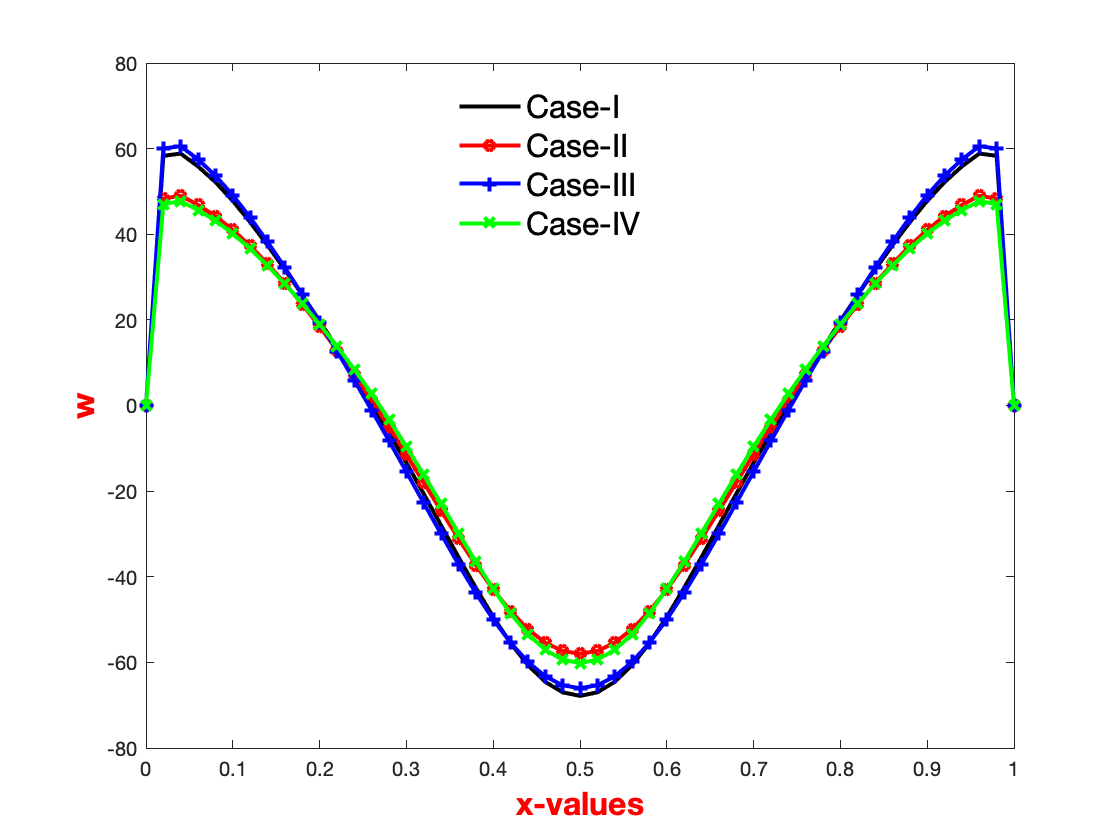}
    \caption{$w$ along horizontal midline} 
\end{subfigure}
\quad 
\begin{subfigure}{0.4\textwidth}
    \includegraphics[width=\textwidth]{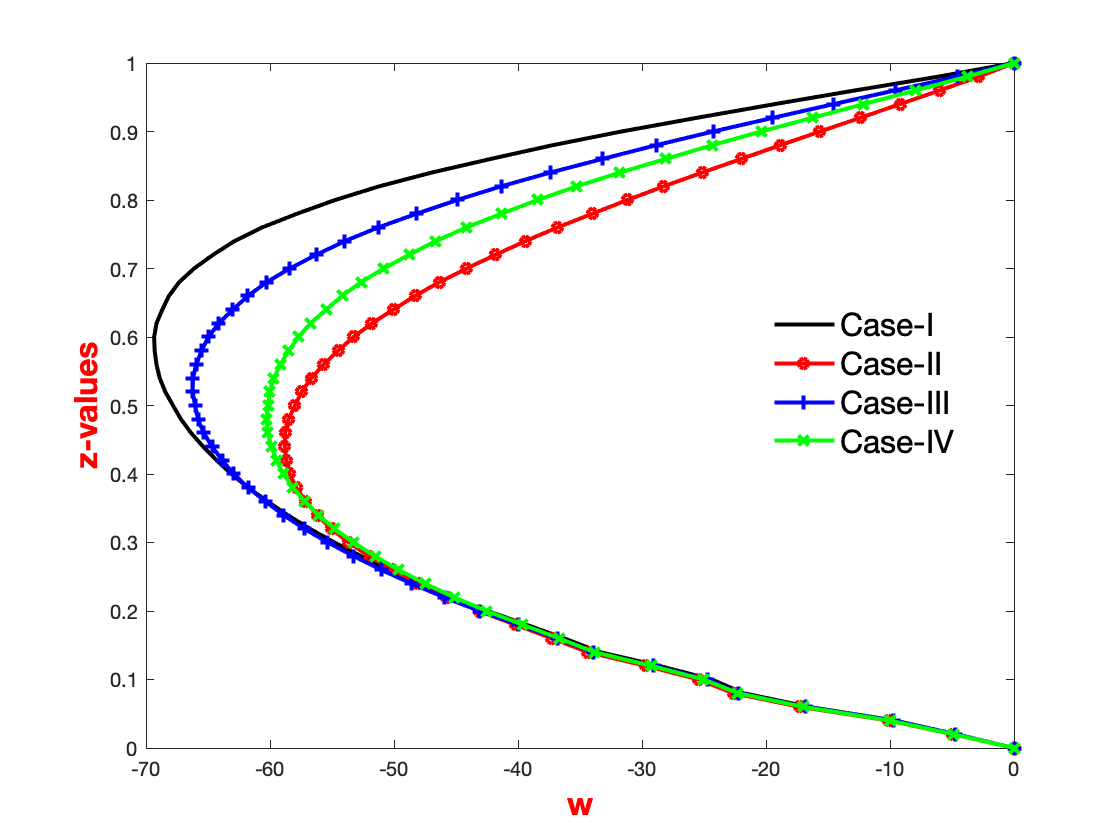}
    \caption{$w$ along vertical midline}
\end{subfigure}
\caption{Velocity components along the horizontal and vertical midlines for all the cases and fixed Rayleigh number $Ra=500$.  }
\label{fig_midline_all_cases}
\end{figure}

Figure~\ref{fig_midline_all_cases} depicts the velocity components along the horizontal and vertical midline in the domain for a fixed value of the Rayleigh number. It is clear from these figures that the magnitude of the velocity is largest along the right and left boundary and symmetric along the vertical midline, thereby creating two convection cells. he computations are done for the fixed Rayleigh number. For all the cases, the vertical velocity $w$ is symmetric about the vertical line $x = 0.5$ and two maximums near the vertical boundary, i.e. $x=0$ and $x=1$, which indicates that for $Ra=500$ there will be two convective cells and the heat transfer rate is maximum at the vertical midline.

\begin{figure}[H]
\centering
\begin{subfigure}{0.4\textwidth}
    \includegraphics[width=\textwidth]{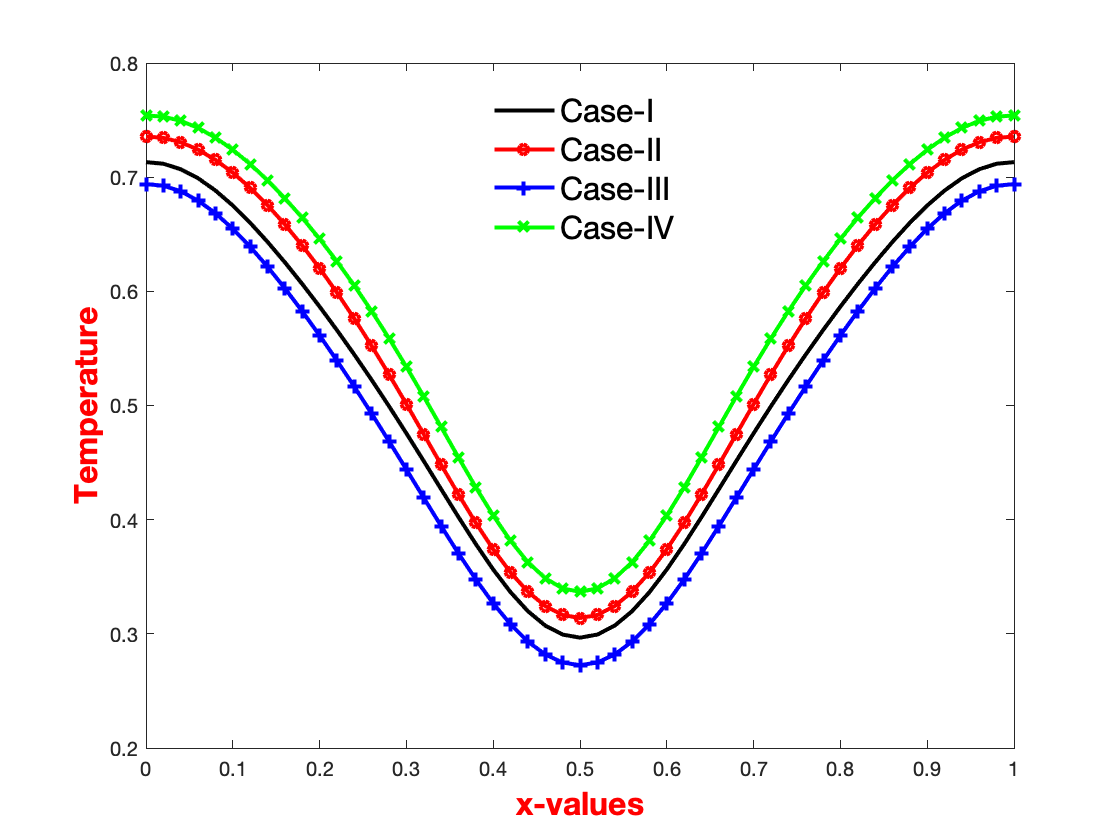}
    \caption{Temperature on the horizontal midline }
\end{subfigure}
\quad 
\begin{subfigure}{0.4\textwidth}
    \includegraphics[width=\textwidth]{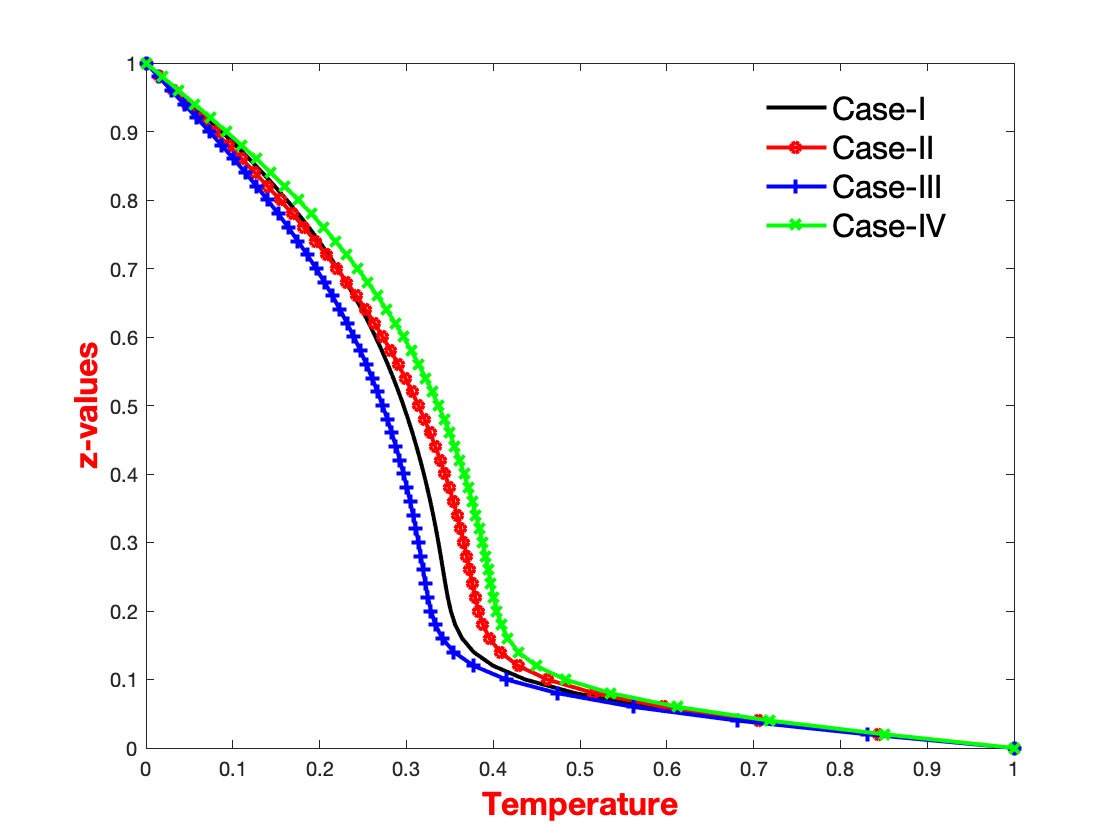}
    \caption{Temperature on the vertical midline} 
\end{subfigure}
\caption{Temperature on the midlines in the domain for a fixed Rayleigh number $Ra=500$.  }
\label{fig_temp_all_cases}
\end{figure}

Figure~\ref{fig_temp_all_cases} presents the temperature values along the midlines in the domain. The temperature profiles show a bigger variation for the quadratic hydraulic resistivity, however the opposite behavior for the quadratic thermal diffusivity.

\section{Conclusion}
The main objective of the present work is to formulate a stable finite element method for the discretization of the Darcy-B\'enard convection model and characterize the effects of hydraulic resistivity and thermal conductivity. The continuity-momentum system of partial differential equations is discretized by using \textit{inf-sup} stable pair of classical  \textit{Taylor-Hood} elements to approximate the velocity-pressure variables. Furthermore, we use a decoupling technique to couple with the time-dependent temperature equation. The overall algorithm is shown to preserve optimal order convergence rates for a ``manufactured solution''.  The numerical computations are assessed for accuracy by comparing the average Nusselt numbers with those obtained from the finite difference method (FDM). Results for the average Nusselt number are in good agreement with the previous results obtained of natural convection in Darcy porous medium using FDM. Simulations are also done for various computational meshes and obtain convergent values of the Nusselt number. 

A key objective of this investigation is to perform finite element simulations of hydro-thermal convection problems. 
 Extensive simulations are performed to study the influence of various cases of hydraulic resistivity and thermal diffusivity. The solution features of flow and temperature characteristics are depicted through velocity vectors and isotherms. A key finding from current research is that the heat transfer rate is the largest for linear variation in hydraulic resistivity and quadratic thermal diffusivity. It is also observed that after the onset of convection, the heat line contours get compressed towards adiabatic vertical walls, deformed vertically, and flattened towards the bottom hot boundary. The fluid circulation is observed to form two convective cells for each Rayleigh number and varying porous medium parameters. 
 
The current computational architecture, the finite element-based numerical method, and the in-house computational code can be easily extended to three-dimensional problems, B\'enard convection in other porous mediums, thermal convection in the nano-fluidic environment, etc. The current numerical method can be extended to include local mesh adaptivity by finding suitable \textit{a-priori} or residual-based refinement indicator. All the aforementioned problems are interesting future works.

\bibliographystyle{plain}
\bibliography{fluids}

\end{document}